\documentclass[graybox]{svmult}
\usepackage{mathptmx}       
\usepackage{helvet}         
\usepackage{courier}        
\usepackage{makeidx}         
\usepackage{graphicx}        
\usepackage{multicol}        
\usepackage[bottom]{footmisc}
\usepackage{amsmath,amssymb,amsfonts}        
\makeindex             
\usepackage{enumitem}
\usepackage{tikz}
\usepackage{tikz-cd}
\usepackage{xspace}
\usepackage{dsfont}
\usepackage{bm}
\usepackage{color}
\usepackage{subfigure}
\usepackage{array}
\usepackage{float}
\usepackage{calc}
\usepackage{import}

\newcommand{\cA}{{\mathcal A}}
\newcommand{\cX}{{\mathcal X}}

\newcommand{\fd}{\mathfrak{d}}
\newcommand{\ff}{\mathfrak{f}}

\newcommand{\fD}{\mathfrak{D}}

\newcommand{\seed}{\textbf{s}}
\newcommand{\Aprin}{\cA_{\text{prin}}}

\newcommand{\lrp}[1]{\left(#1\right)}

\newcommand{\C}{\mathbb{C}}
\newcommand{\R}{\mathbb{R}}
\newcommand{\Q}{\mathbb{Q}}
\newcommand{\Z}{\mathbb{Z}}

\newcommand{\rank}{\operatorname{rank}}
\DeclareMathOperator{\Hom}{Hom}
\newcommand{\Spec}{\operatorname{Spec}}
\newcommand{\can}{\operatorname{can}}

\newcommand{\mH}{\mathcal{H}}

\newcommand{\CP}{\mathbb{C}P}

\newcommand{\PxP}{\mathbb{C}P^1 \times \mathbb{C}P^1}

\newcommand{\BlI}{\mathbb{C}P^2\# \overline{\mathbb{C}P^2}}

\newcommand{\BlIII}{\mathbb{C}P^2\# 3\overline{\mathbb{C}P^2}}

\newcommand{\BlV}{\mathbb{C}P^2\# 5\overline{\mathbb{C}P^2}}
\newcommand{\BlVI}{\mathbb{C}P^2\# 6\overline{\mathbb{C}P^2}}

\newcommand{\val}{\mathrm{val}}

\newcommand{\del}{\partial}

\newcommand{\ind}{\mathrm{ind}}

\newcommand{\tY}{\tilde{Y}}

\newcommand{\bY}{\overline{Y}}

\newcommand{\Proj}{\mathrm{Proj}}
\newcommand{\TV}{\mbox{TV}}

\usepackage{upgreek}
\usepackage{caption}
\usetikzlibrary{shapes,arrows,patterns,decorations.pathreplacing,shapes.geometric}
\tikzstyle{sing}=[star, draw=black, fill=white, inner sep=0pt, minimum size=5pt]

\begin{document}
\title*{Algebraic and symplectic viewpoint on compactifications of two-dimensional cluster varieties of finite type}
\titlerunning{Compactifications of cluster varieties}
\author{Man-Wai Mandy Cheung and Renato Vianna}
\institute{Man-Wai Mandy Cheung \at Harvard University, One Oxford Street, Cambridge, MA 02138, United States of America, \\ \email{mwcheung@math.harvard.edu}
\and Renato Vianna \at Institute of Mathematics, Federal University of Rio de Janeiro; Av. Athos da Silveira Ramos, 149 - Ilha do Fundão, Rio de Janeiro - RJ, 21941-909, Brazil, \\ \email{renato@im.ufrj.br}}
%
%
\maketitle

\abstract{In this article we explore compactifications of cluster varieties of finite type in complex dimension two. Cluster varieties can be viewed as the spec of a ring generated by theta functions and a compactification of such varieties can be given by a grading on that ring, which can be described by positive polytopes \cite{ghkk}. 
In the examples we explore, the cluster variety can be interpreted as the complement of certain divisors in del Pezzo surfaces. 
In the symplectic viewpoint, they can be described via almost toric fibrations over $\R^2$ (after completion). Once identifying 
them as almost toric manifolds, one can symplectically view them inside other del Pezzo surfaces. So we can identify other symplectic
compactifications of the same cluster variety, which we expect should also correspond to different algebraic compactifications. 
Both viewpoints are presented here and several compactifications have their corresponding polytopes compared. The finiteness of the
cluster mutations are explored to provide cycles in the graph describing monotone Lagrangian tori in del Pezzo surfaces connected via 
almost toric mutation \cite{Vi16a}. 
}

\section{Introduction}

Cluster algebras, introduced by Fomin and Zelevinsky \cite{cluster1}, are subalgebras of rational functions in $n$ variables.
The generators of cluster algebras are called the cluster variables. 
Instead of being given the complete sets of generators and relations as other commutative rings, a cluster algebra is defined from an (initial) seed, which includes a set of the generators and a matrix.  
An iterative procedure called mutation would produce new seeds from a given seed and this process gives all the cluster variables.
The cluster algebra is then defined to be the ring generated by all cluster variables. 

Geometrically, the cluster varieties, described by Fock and Goncharov \cite{FGdilog1}, and by Gross, Hacking, Keel in \cite{GHK_bir}, are defined in a similar manner.
A seed data now would be associated to an algebraic torus.
The mutation procedures give the birational transformations used to glue the tori. 
A cluster variety is then the union of the tori under the gluing.

The compactification of the cluster varieties can be given by a Rees construction.
Combinatorially, the construction can be described by `convex' polytopes, called the positive polytopes \cite{ghkk}. 
The article \cite{cpt} showed that the positive polytopes satisfy a convexity condition called `broken line convexity'. 
As the seed mutates, the polytope mutates correspondingly. 
One can then give a mutation process to the polytopes. 
Note that under this type of mutation, there is no change in the compactification. 

More generally, one can similarly describe the compactification of the log Calabi-Yau surfaces studied in \cite{GHKlog}. 
In this case, one would construct the dual intersection complex of a given Looijenga pair. 
The underlying topological space of the complex will carry an affine manifold structure. 
The affine structures would correspond to another type of mutation for the positive polytopes. 

On the other hand, in the symplectic viewpoint, mutations were exploited in four dimensional symplectic geometry \cite{Vi14,Vi16a}, inspired by the pioneering work of
Galkin-Usnhish \cite{GaUs10} (further developed in \cite{ACGK12}), and being grounded on the development of almost toric fibrations (ATFs)
by Symington \cite{Sy03}. Upon identifying an almost toric fibration of a open variety, we can symplectically identify it as 
a symplectic submanifold of some closed symplectic manifold. We will refer to it as a (symplectic) compactification. In the examples
of this paper, we can identify the symplectic form as the K\"{a}hler form of del Pezzo surfaces. We expect that symplectic compactifications
can be translated to algebraic compactifications under certain nuances discussed in Section \ref{sec:symp}.

This paper is an attempt to understand the two notions.
The motivation of both sides come from the Strominger-Yau-Zaslow conjecture -- the conjecture suggests there are special Lagrangian fibrations for the Calabi-Yau manifold and its mirror space over the base $B$. 
The construction of the log-CY variety from the symplectic side is via the almost toric fibration, 
Meanwhile, in the algebro-geometric side, the construction can be described in terms of the wall crossing structures called the scattering diagrams. 

We begin with the algebro-geometric perspective in Section \ref{sec:ag}. 
In this section, we will discuss cluster varieties, positive polytopes, compactifications, and the mutations of the polytopes.
Then in Section \ref{sec:SYZ}, we give a perspective on how cluster varieties
and scattering diagrams arise from considering wall crossing corrections 
as one attempt to build a mirror in terms of the SYZ picture. In complex dimension two, the wall-crossing happens when we consider singular Lagrangian fibrations known as almost toric fibrations (ATFs).
In particular, we illustrate the idea in terms of the $A_2$ cluster variety -- compactified as the del Pezzo surface of degree 5 in Section \ref{sec:A2}.
Afterward, we explore compactifications of cluster varieties using the almost toric viewpoint in Section \ref{sec:symp_comp}. 

The symplectic geometry approach to compactification via almost toric fibrations makes no reference to the complex structure, while the algebro-geometric approach does not fix a symplectic form. Nonetheless, because a similar set of data can encode the scattering
diagram as well as an ATF, we seem to always be able to relate compactifications, encoded by the same polytope in both pictures. 
We aim to show the correspondence between the symplectic compactification of the cluster varieties to the algebro-geometric version in our upcoming papers.

\begin{acknowledgement}
This project initialized from discussions during the conference "Tropical Geometry And Mirror Symmetry" in the MATRIX Institute. 
The authors would like to thank the  MATRIX institute for their hospitality.

The authors would like to thank Denis Auroux, and the referee for helpful feedback on the first version of the paper. 
The first author would like to thank Tim Magee, and Yu-shen Lin for helpful discussions. The first author is supported by NSF grant DMS-1854512.
The second author is supported by Brazil's National Council
of scientific and technological development CNPq, via the research fellowships
405379/2018-8 and 306439/2018-2, and by the Serrapilheira Institute grant
Serra-R-1811-25965.

\end{acknowledgement}

\section{Mutations in algebraic geometry} \label{sec:ag}

\subsection{Cluster varieties} \label{sec:ClusterVar}

We will first recall some notation used in the definition of a cluster varieties. 
A \emph{fixed data} consists of a lattice $N$ with a skew-symmetric bilinear form $\{ \cdot , \cdot \} : N \times N \rightarrow \Q$,
an index set $I$ with $|I| = \rank N$,
positive integers $d_i$ for $i \in I$,
a sublattice $N^{\circ} \subseteq N$ of finite index with some integral properties,  the dual lattice $M = \Hom (N, \Z)$ and the corresponding $M^{\circ} = \Hom(N^{\circ}, \Z)$. 
One can refer to \cite{GHK_bir} for the full definition of fixed data. 
Consider $N_{\R} = N \otimes \R$ and $M_{\R} = M \otimes \R$. 

Given this fixed data, a $\emph{seed data}$ for this fixed data is $\seed := ( e_i \in N \mid i \in I)$,
where  $\{ e_i \}$ is a basis for $N$.
The basis for $M^{\circ}$ would then be $f_i = \frac{1}{d_i} e_i^*$ . 
One can then associate the seed tori 
\[\cA_{\seed} = T_{N^{\circ}} = \Spec \Bbbk [M^{\circ}], \quad
\quad \cX_{\seed} = T_{M} = \Spec \Bbbk [N].\]
We will denote the coordinates as $X_i = z^{e_i}$ and $A_i = z^{f_i}$ and they are called the \emph{cluster variables}. 
Similar to the definition of cluster algebras, there is a procedure, called \emph{mutation}, to produce a new seed data $\mu(\seed)$ from a given seed $\seed$. 
The mutation formula is stated in \cite[Equation 2.3]{GHK_bir} which we will skip here. 
The essence is that we will obtain new seed tori $\cA_{\mu(\seed)}$, $\cX_{\mu(\seed)}$ from the mutated seed. 
Between the tori, there are birational maps $\mu_{\cX}: \cX_{\seed} \dashrightarrow \cX_{\mu(\seed)}$, $\mu_{\cA}: \cA_{\seed} \dashrightarrow \cA_{\mu(\seed)}$ which are stated in \cite[Equations 2.5, 2.6]{GHK_bir}.
Note that those birational maps are basically the mutations of cluster variables as in Fomin and Zelevinsky \cite{cluster1}. 

Let $\cA$ be an union of tori glued by $\cA$-mutation $\mu_{\cA}$.
A smooth scheme $V$ is a {\emph{cluster variety of type $\cA$}} if there is a birational map $ \mu: V \dashrightarrow \cA$ which is an isomorphism outside codimension two subsets of the domain and range.
The {\emph{cluster variety of type $\cX$}} is defined analogously. 

The $\cA$ and $\cX$ cluster varieties can be fit into the formalism of the cluster varieties with principal coefficients $\Aprin$. The scheme $\Aprin$ is defined similarly to the $\cA$ by `doubling' the fixed data, i.e. considering $\widetilde{N} = N \oplus M^{\circ}$ as fixed data as in \cite[Construction 2.11]{GHK_bir}. 
Then there are two natural inclusions. The first one is 
\begin{align*}
    \widetilde{p}^* \colon  &N \rightarrow \widetilde{M}^{\circ} = M^{\circ} \oplus N,  \\ 
    &n \mapsto (p^*(n),n),
\end{align*}
where $p^*(n) = \{ n, \cdot \} \in M^{\circ}$ in the case of no frozen variable. 
Then for any seed $\seed$, 
note that $\cA_{\text{prin}, \seed} = T_{\widetilde{N}^{\circ}} $, and $ \cX_{\seed} = T_M$, the there is the exact sequence of tori
\[
1 \rightarrow T_{N^{\circ}} \rightarrow \cA_{\text{prin}, \seed} \xrightarrow{\widetilde{p}}  \cX_{\seed} \rightarrow 1. 
\]
The map $\widetilde{p}$ commutes with the mutation maps and thus we get the morphism $\widetilde{p} \colon \Aprin \rightarrow \cX$. 
Further the $T_{N^{\circ}}$ action on $\cA_{\text{prin}, \seed}$ extends to $\Aprin$ which makes $\widetilde{p} $ a quotient map. 
Thus, the $\cX$ variety can be seen as $\Aprin / T_{N^{\circ}}$.

The second inclusion is 
\begin{align*}
    \pi^* \colon &N \rightarrow M^{\circ}, \\
    &n \mapsto (0,n). 
\end{align*}
In this case, the $\pi^*$ map induces a projection $\pi: \Aprin \rightarrow T_M$. 
Then the usual $\cA$ variety is $\pi^{-1}(e)$, where $e$ is the identity of $T_M$. 

We would like to indicate another viewpoint of the cluster varieties here. 
The mutation maps may be described in terms of elementary transformation of $\mathbb{P}^1$ bundles. 
Thus the cluster varieties can also be seen as the blowups of toric varieties (up to codimension two) as well. 

Given a seed data, consider the fans 
\[
\Sigma_{\seed, \cA} := \{ 0 \} \cup \{ \R_{\geq 0} d_i e_i \mid i \in I \} \subseteq N^{\circ}, \
\Sigma_{\seed, \cX} := \{ 0 \} \cup \{- \R_{\geq 0} d_i v_i \mid i \in I \} \subseteq M,
\]
where $v_i = p^*(e_i)$ and the $i$ only runs over the unfrozen variables if the frozen variables exist. 
Let $\TV_{\seed, \cA}$ and $\TV_{\seed, \cX}$ be the respective toric varieties. 
Denote $D_i$ to be the toric divisor corresponding to the one-dimensional ray in one of these fans. 
Define the closed subschemes
\[
Z_{\cA, i} := D_i \cap \bar{V} (1+z^{v_i}) \subseteq \Sigma_{\seed, \cA}, \
Z_{\cX, i} := D_i \cap \bar{V} \left((1+z^{e_i})^{\ind \ d_i v_i}\right) \subseteq \Sigma_{\seed, \cX},
\]
where $\bar{V}$ denote the closure of the variety $V$, and $\ind  \ d_i v_i$ is the greatest degree of divisibility of $  d_i v_i$ in $M$. 
Then consider the pairs $(\widetilde{\TV}_{\seed, \cA}, D)$ and $(\widetilde{\TV}_{\seed, \cX}, D)$ consisting of the blowups of $\TV_{\seed, \cA}$ and $\TV_{\seed, \cX}$ respectively, with $D$ the proper transform of the toric boundaries. 
Define $X_{\seed, \cA} = \widetilde{\TV}_{\seed, \cA} \setminus D$
and $X_{\seed, \cX} = \widetilde{\TV}_{\seed, \cX} \setminus D$. 
When the seed $\seed$ mutates to $\seed'$, the corresponding $X_{\seed, \cA}$, $X_{\seed', \cA}$ and $X_{\seed, \cX}$, $X_{\seed', \cX}$ are isomorphic outside a codimension two set. 
In finite type, where there are only finitely cluster variables,  the $\cA$ and $\cX$ would then also isomorphic to $X_{\seed, \cA}$ and $X_{\seed, \cX}$. 
Note that the whole set up here is building a toric model for the cluster varieties. 
We will introduce the notion of toric model for log Calabi Yau surfaces later in Section \ref{sec:can}.

\subsubsection*{Scattering diagrams}

Scattering diagrams live in the tropicalization of the cluster varieties. 
One can also see the diagrams encode the structure of the cluster varieties combinatorially. 

A {\it{wall}} in $M_{\R}$ is a pair $(\fd, f_{\fd})$ where
 $\fd \subseteq M_{\R}$ is a convex rational polyhedral cone of codimension one, contained in $n^{\perp}$ for some $n \in N$, and 
 $f_{\fd}  = 1+ \sum_{k \geq 1} c_k z^{kp^*(n)}$, where $c_k \in \C$. 
A wall $(\fd, f_{\fd})$ is called \emph{incoming} if $p^* (n) \in \fd$. Otherwise it is called \emph{outgoing}.
A \emph{scattering diagram} $\fD$ is then a collection of walls with certain finiteness properties. 
Given a seed, an $\Aprin$-cluster scattering diagram can be constructed \cite{ghkk} and canonically determined by this given seed data.
The $\cA$ scattering diagram can be obtained by the projection $\widetilde{M}_{\R} \rightarrow M_{\R}$ while 
the $\cX$ scattering diagrams can be defined as slicing the $\Aprin$ scattering diagrams by considering $\{ (m,n) \mid m = p^*(n) \}$. 

It is worth addressing here that for finite type, each chamber, i.e. the maximal cone, of the scattering diagram can be associated to a torus. The wall functions $f_\fd$ are actually representing the birational maps between the tori. 
Thus the cluster varieties can be seen as gluing of tori associated to the chamber via the wall crossing. 

\begin{figure}[h!]   
  \begin{center} 
 \centerline{\includegraphics[scale= .7]{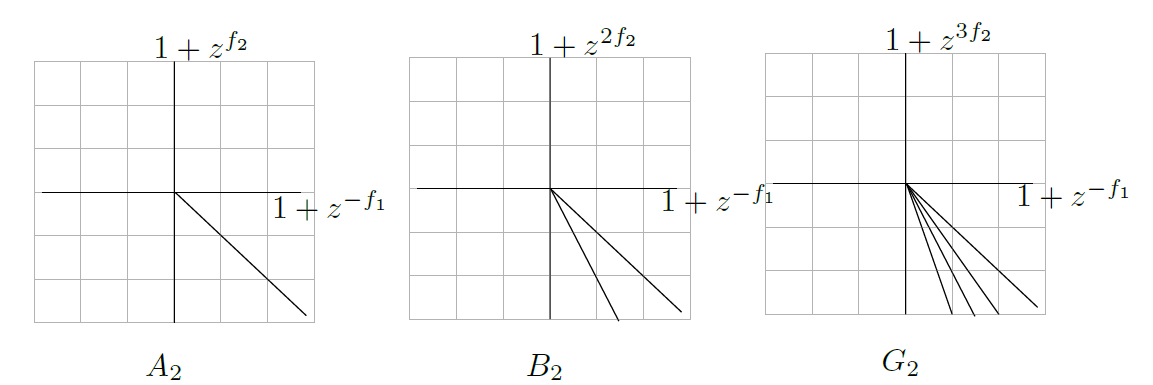}}
\caption{$\cA$-scattering diagrams for rank 2 finite type.} \label{fig:rank2}
\end{center}
\end{figure}

\begin{figure}[h!]   
  \begin{center} 
 \centerline{\includegraphics[scale= .7]{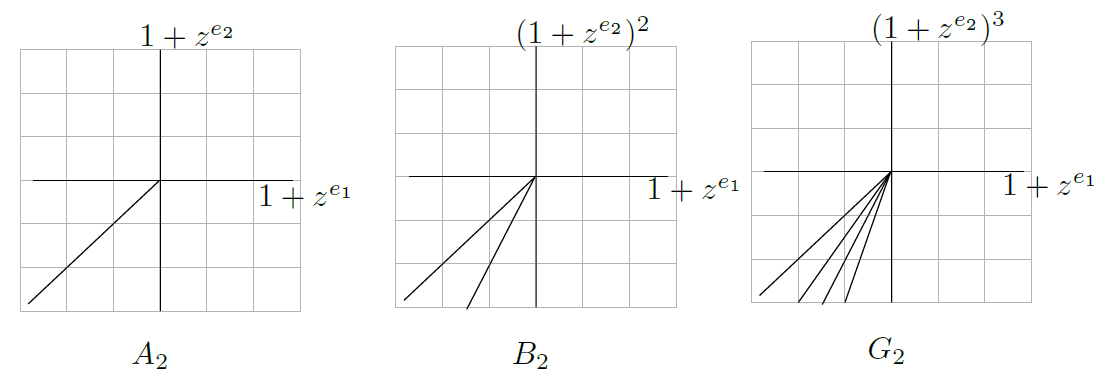}}
\caption{$\cX$-scattering diagrams for rank 2 finite type.} \label{fig:Xrank2}
\end{center}
\end{figure}

In this article, we will focus on the dimension 2 cluster varieties of finite type.
The $\cA$ scattering diagrams are listed as in Figure \ref{fig:rank2} while the $\cX$ scattering diagrams of rank 2 finite type are listed as in Figure \ref{fig:Xrank2}.

\subsubsection*{Mutation of scattering diagrams}

As noted in the previous section, a seed determines canonically a scattering diagram. 
Two mutation-equivalent seeds would then give two different scattering diagrams. 
It is natural to consider `mutation equivalent' scattering diagrams. 
This equivalence is given by piecewise linear maps on the lattices which are very similar to those in Section \ref{sec:symp} and hence we will state here. 

Consider two seeds $\seed$ and $\seed'$ which are just one mutation step apart, i.e. $\seed' = \mu_k (\seed)$ for some $k \in I$.
Then the corresponding scattering diagrams $\fD_{\seed}$ and $\fD_{\seed'}$ are equivalent to each other by the transformation $T_k: M^{\circ} \rightarrow M^{\circ}$,
\begin{equation} \label{eqn:semilinear}
T_k(m) = \left\{
\begin{array}{l r}
m + \langle d_k e_k , m \rangle v_k,   & \text{ for }m \in \mH_{k, +}\\
m, & \text{ for } m \in \mH_{k,-}
\end{array}
\right.
\end{equation}
for $m \in M^{\circ}$, $v_k = p^*(e_k)$, and
$\mH_{k, +} = \{ m \in M_{\R} | \langle e_k , m \rangle \geq 0 \}$,  $
	\mH_{k, -} = \{ m \in M_{\R} | \langle e_k , m \rangle \leq 0 \}$.
Extending $T_k$ to the wall functions  \cite[Theorem 1.24]{ghkk} will lead us to another consistent scattering diagram $T_k (\fD_s)$ which is shown to be equivalent to $\fD_{\mu_k(s)}$.

\begin{figure}[H]   
  \begin{center} 
 \centerline{\includegraphics[scale= 0.3]{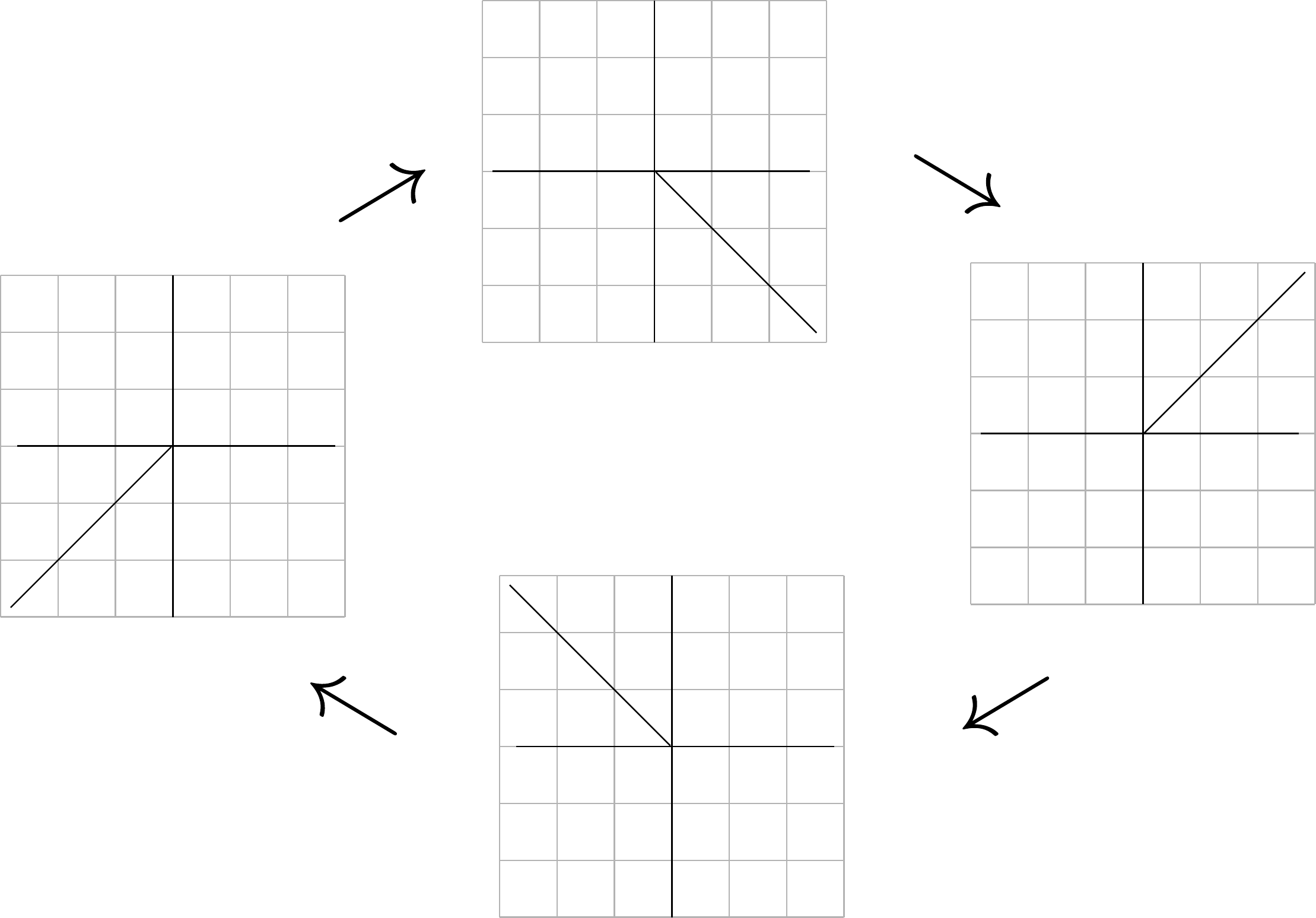}}
\caption{Mutation of the $\cX$ scattering diagram of type $A_2$ starting at the index $1$} \label{fig:XA2mutate}
\end{center} 
\end{figure}

\begin{figure}[h!]   
  \begin{center} 
 \centerline{\includegraphics[scale=2]{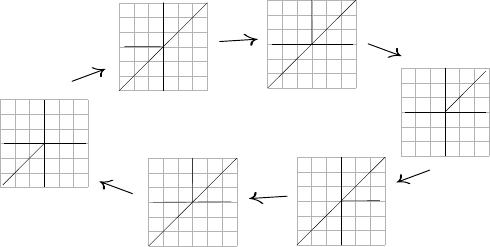}}
\caption{Mutation of the $\cX$ scattering diagram of type $A_2$ starting at the index $2$}  \label{fig:A2mutate_2}
\end{center}
\end{figure}

We can similarly define the mutation for the $\cX$ scattering diagrams from $\Aprin$.
For the scattering diagram of type $A_2$ in Figure \ref{fig:Xrank2}, we can obtain the mutation process for the $\cX$ scattering diagram as in Figure \ref{fig:XA2mutate} and Figure \ref{fig:A2mutate_2}.  

Note that the scattering diagrams are determined by seeds while the mutation of seeds are given by blow ups and blow downs of toric varieties. 
Thus the mutation of scattering diagrams actually represents this procedure of blowups and blowdowns. We are going to discuss a similar construction in Section \ref{sec:symp} with a symplectic perspective.

\subsubsection*{Theta functions and the canonical algebras}

Theta functions give the generators of the canonical basis of the cluster algebras. 
Given the cluster variety $V = \cA, \Aprin, \cX$, the corresponding character lattice is $L = M^{\circ}$, $\widetilde{M}^{\circ}$, or $N$. 
A theta function $\vartheta_p$ is associated to each point $p \in L$ by a combinatorial object -- broken lines which are piecewise linear paths in $L_{\R}$ together with decorating monomials at each linear segment. 

The free module generated by theta functions is endowed with an algebra structure from the multiplication between theta functions.
Indeed the structure constants in the multiplications of theta functions can be given in terms of counting broken lines.
The product of two theta functions can be expressed as
\begin{align} \label{eq:product}
    \vartheta_p \cdot \vartheta_q = \sum_r \alpha(p,q,r) \vartheta_r,
\end{align}
where the structure constants $\alpha\lrp{p,q,r}$  can be explicitly defined by counting broken lines with certain boundary conditions \cite[Proposition 6.4]{ghkk}.
In this finite type case, the structure constants $\alpha$ define (\cite[Corollary 8.18]{ghkk}) the finitely generated $\C$-algebra structure on 
\[
\can (V) := \bigoplus_{r \in L } \C \cdot \vartheta_r.
\]
We will then define $X:= \Spec (\can (V))$.

\subsection{Positive polytopes}

With the multiplication structure of the theta functions, we can now state the definition of a positive set-- the property required for a set and its dilations to define a graded ring.

For $S \subseteq L_{\R} = L \otimes \R$ a closed subset, define the cone of $S$ as
\[
\textbf{C}(S) = \overline{ \{ (p,r) \mid p \in rS, r \in \R_{\geq 0} \}} \subset  L_{\R} \times \R_{\geq 0}. 
\]
Denote $d S(\Z) = \mathbf{C} (S) \cap \left( L \times \{d\}\right)$ which is viewed as a subset of $L$. 

A closed subset $S \subset L_{\R}$ is called \emph{positive} if for any non-negative integers $d_1$, $d_2$, any $p_1 \in d_1 S(\Z)$, $p_2 \in d_2S(\Z)$, and any $r \in L$ with $\alpha (p_1 ,p_2, r) \neq 0$, then $r \in (d_1 + d_2 ) S(\Z)$.

\begin{figure}  
  \begin{center} 
 \centerline{\includegraphics[scale= 0.45]{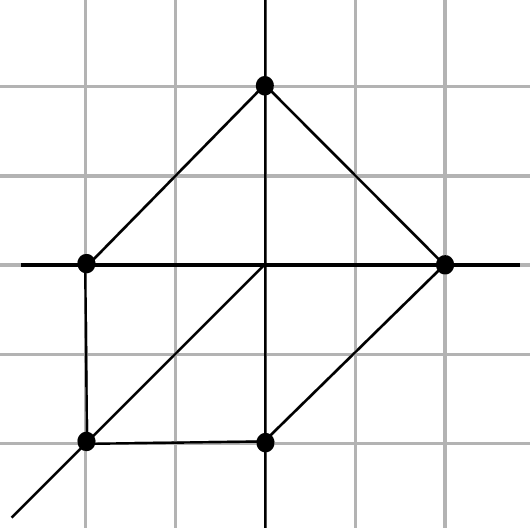}}
 \caption{Positive polytope of $\cX$ cluster variety of type $A_2$.} \label{fig:A2polytope}
\end{center} 
\end{figure}

In the ongoing example of cluster varieties of type $A_2$, we consider the polytope with vertices $(1,0), (0,1), (-1,0), (0,-1), (1,-1)$ as indicated in Figure \ref{fig:A2polytope}.
Note that this polytope is in the $\cX$ diagram thus there is a flip from Figure \ref{fig:dp5_00}.
This polytope is indeed positive. 
In Section \ref{sec:A2}, there is a detail discussion of such a polytope in the $\cA$ side. 
A similar calculation in this $\cX$ case will still hold, thus this will correspond to the del Pezzo surface of degree 5 \cite{ghkk}. 

\begin{figure}[h!]   
  \begin{center} 
 \centerline{\includegraphics[scale= 0.4]{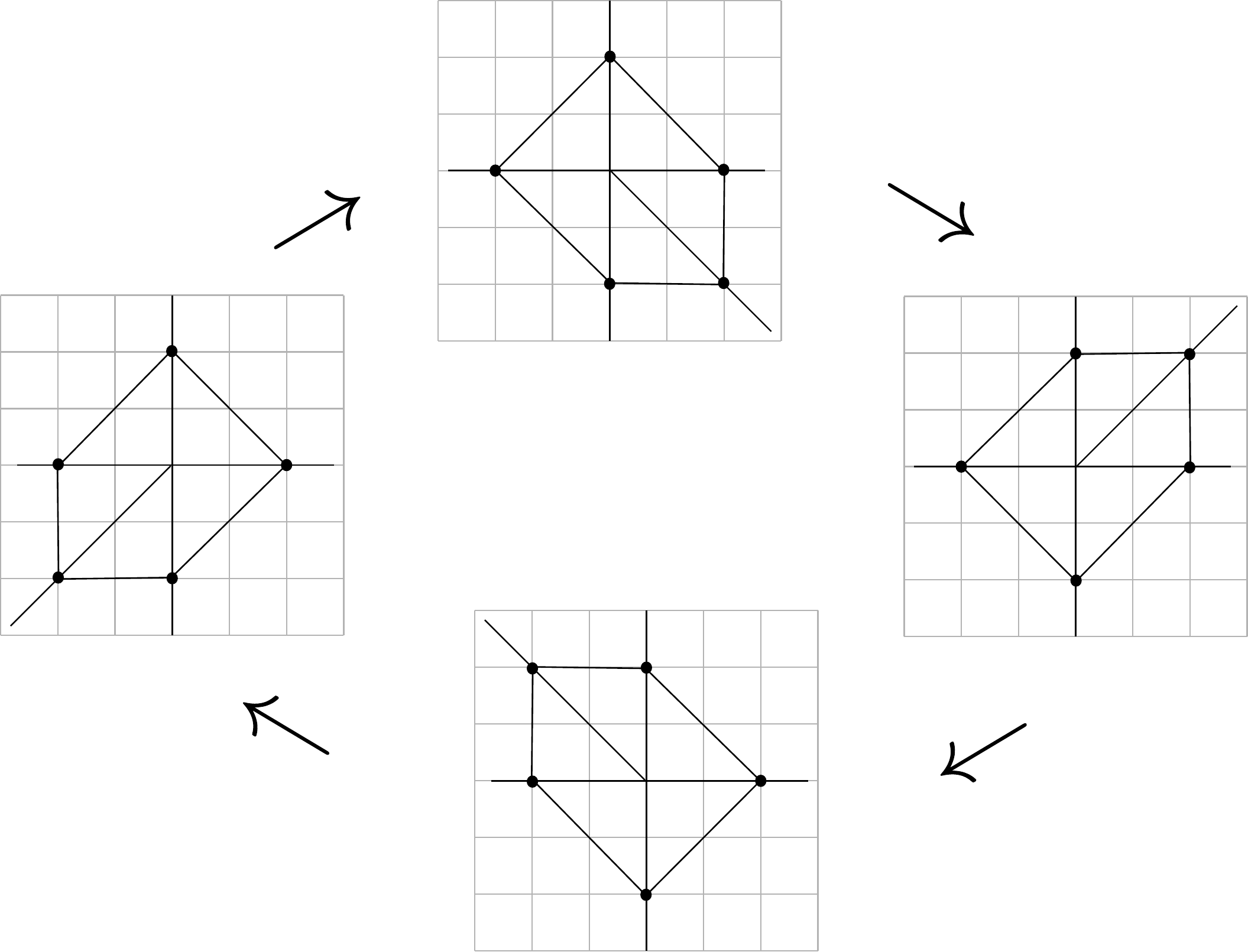}}
\caption{Mutation of polytope starting from index $1$}   \label{fig:Xpolymutate_1}
\end{center} 
\end{figure}

\begin{figure}[h!]   
  \begin{center} 
 \centerline{\includegraphics[scale= 2]{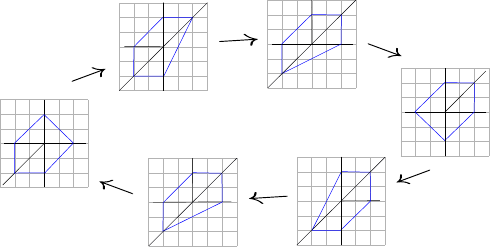}}
\caption{Mutation of polytope starting from index $2$} 
\label{fig:whatever}
\end{center} 
\end{figure}

We can apply the mutation sequences in Figure 
\ref{fig:XA2mutate} and \ref{fig:A2mutate_2} to the polytope in Figure \ref{fig:A2polytope}. 
Mutations of the polytopes as in Figure \ref{fig:Xpolymutate_1} and \ref{fig:whatever} will be obtained respectively. 
In the next section, we will describe mutations of the polytopes from a symplectic point of view. 
We observe that the mutation sequences of polytope in Figure \ref{fig:Xpolymutate_1} and \ref{fig:whatever} are the same as the sequences in Figure \ref{fig:dp5_2_1} and \ref{fig:dp5_2_2} respectively. 
The cluster mutation of the scattering diagrams comes from a change of seed data, i.e. a change of the initial variables. 
Thus the underlying spaces are all isomorphic. 

\subsubsection*{Compactifications from positive polytopes}

We will roughly go over the geometric meaning behind the positive polytopes in this section. 
The motivation can be seen as the construction of projective toric varieties from the convex polytopes. 

For rank 2 cluster varieties, since the $\cA$ scattering diagrams are well defined, we can consider $\bar{S}$ the positive polytopes in the $\cA$ scattering diagrams. 
For this set $\bar{S}$, define $\widetilde{S} = \bar{S} + N_{\R}$ which is obviously positive. 
Thus we can define the graded ring
\begin{align*}
    \widetilde{R}_{\widetilde{S}} = \bigoplus_{d \geq 0} \bigoplus_{q \in d \widetilde{S}(\Z)} \C \vartheta_q x^d \subset \can(\Aprin) [x],
\end{align*}
with grading defined by $x$.

Define $Y_{\Aprin}:=\Proj (\widetilde{R}_S) \rightarrow T_M$.
For the $\cA$ variety, we take $Y_{\cA}$ as the fiber over $e \in T_M$ in this map. 
More generally, for the $\cA_t$ variety, $t \in T_M$, we can still take $Y_{\cA_t}$ as the fiber over $t \in T_M$. 
Consider $X=\Spec (\can (V))$,  for $V = \Aprin, \cA_t$, as in the previous subsection. 
Define $B = Y \setminus X$. 
Then \cite{ghkk} showed that, $X$ is a Gorenstein scheme with trivial dualizing sheaf, in particular, for $V = \Aprin, \cA_t$, $X$ is a $K$-trivial Gorenstein log canonical variety. 
In this finite rank 2 case, for $V = \Aprin, \cA$, $X \subseteq Y$ is a minimal model, i.e. $Y$ is a projective normal variety, $B \subset Y$ is a reduced Weil divisor, $K_Y + B$ is trivial, and $(Y, B)$ is log canonical. 

For the case of the $\cX$ varieties, as indicated in Section \ref{sec:ClusterVar}, the $\cX$ varieties are quotients of the $\Aprin$ varieties. 
Thus we will consider still consider $\Aprin$ but instead see the lattice as $\widetilde{M}^{\circ}$ instead of $\widetilde{N}$ (which are actually isomorphic). 
We can repeat the same procedure as before and then obtain the compactification of $X= \Spec (\can (\cX))$. 
The scheme $X$ is also a $K$-trivial Gorenstein log canonical variety.

\subsection{Canonical scattering diagrams} \label{sec:can}

In the last section, we note that the cluster mutations of the scattering diagrams are not changing the underlying schemes. 
We are proposing another type of mutation which is given by the monodromy on $B$. 
We are going to understand the ideas behind from the mirror construction suggested by Gross, Hacking, and Keel in \cite{GHKlog}. 
In Section \ref{sec:SYZ}, we will discuss the affine structure and monodromy from the SYZ perspective.

Consider a pair $(Y,D)$, where $Y$ is a smooth rational projective surface, and
$D$ is an anti-canonical cycle of projective lines.
We will call such a pair a Looijenga pair. 
Let $X= Y \setminus D$.
The {tropicalization} of $(Y,D)$ is a pair $(B, \Sigma)$, where $B$ is an integral linear manifold with singularities, and $\Sigma$ is a decomposition of $B$ into cones. 
The pair $(B, \Sigma)$ can be constructed by associating each node $p_{i,i+1} $ of $D$ a rank two lattice with basis $v_i$, $v_{i+1}$. 
Denote the cone generated by $v_i$, $v_{i+1}$ as $\sigma_{i, i+1} \subset M_{i, i+1} \otimes \mathbb{R} $. 
The cones $\sigma_{i, i+1}$ and $\sigma_{i-1, i}$ are  glued over the ray $\rho_i = \mathbb{R}_{\geq 0} v_i$
 to obtain a piecewise linear manifold $B$ homeomorphic to $\mathbb{R}^2$ and $\Sigma = \{ \sigma_{i, i+1} \} \cup \{ \rho_i \} \cup \{0\} $.
 
The integral affine structure on $B_0 = B\setminus \{0\}$ can be defined by the charts 
\[\psi_i: U_i = \mbox{Int}(\sigma_{i-1, i} \cup \sigma_{i, i+1}) \rightarrow M_{\R}, 
\]
where 
\[
\psi_i (v_{i-1}) = (1,0), \ \psi_i (v_{i}) = (0,1), \ \hbox{and} \quad  \psi_i (v_{i+1}) = (-1, -D^2), 
\]
and $\psi_i$ is linear on $\sigma_{i-1, i}$ and $ \sigma_{i, i+1}$.

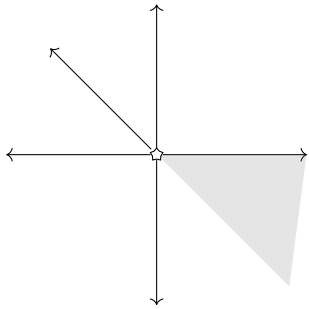
\begin{figure}
\centering
    \begin{tikzpicture}
    \node [inner sep=0] (0) at (0,0) {$\bullet$};
    \draw[->] (0) -- ++(0:2) ; 
    \draw[->] (0)--++(90:2); 
    \draw[->] (0)--++(135:2) ;
    \draw[->] (0)--++(180:2) ; 
    \draw[->] (0)--++(270:2) ; 
     \draw[fill=gray, opacity=0.2, draw=white] (0,0) --(315:2.5)--(360:2) ;
     \fill(0,0) node[sing] {};
\end{tikzpicture}
    \caption{The tropicalization $(B, \Sigma)$ of the del Pezzo surface of degree 5}
    \label{fig:affine}
\end{figure}

Now consider $Y$ the del Pezzo surface of degree 5 and $D$ the anti-canonical cycle of five (-1)-curves. 
The construction of the charts $\psi $ will then give
\[
\psi (v_{1}) = (1, 0), 
\psi (v_{2})  = (0, 1), 
\psi (v_{3})  = (-1, 1),
\psi (v_{4}) = (-1, 0), 
\psi (v_{5})  = (0,-1).
\]
Note however that having $\psi (v_{4}) = (-1, 0), 
\psi (v_{5})  = (0,-1)$ will lead to 
\[\psi(v_1) \leadsto (1,-1), \ \psi(v_2) \leadsto (1,0)\] and this is NOT what we began with: $\psi (v_{1}) = (1, 0)$, and $\psi (v_{2})  = (0, 1)$.
Thus we would like to identify the cone spanned by $(1,0)$ and $(0,1)$, and the cone spanned by $(-1,1)$ and $(1,0)$. 
This introduces the monodromy
\[
(1,0) \mapsto  (1,1), \quad (0,1)\mapsto (1,0),
\]
 to $B_0$. The affine structure is illustrated in Figure \ref{fig:affine}. 

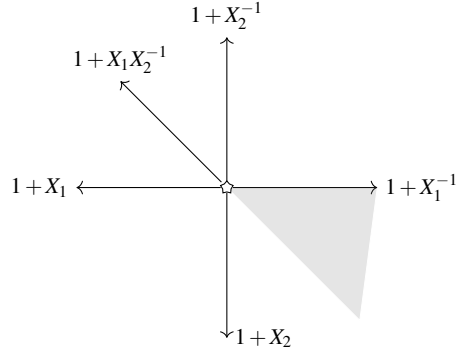
\begin{figure}
\centering
    \begin{tikzpicture}
    \draw[->] (0) -- ++(0:2) node[right]{$1+X_1^{-1}$}; 
    \draw[->] (0)--++(90:2) node[above]{$1+X_2^{-1}$}; 
    \draw[->] (0)--++(135:2) node[above]{$1+X_1X_2^{-1}$};
    \draw[->] (0)--++(180:2)  node[left]{$1+X_1$};  
    \draw[->] (0)--++(270:2)  node[right]{$1+X_2$}; ; 
     \draw[fill=gray, opacity=0.2, draw=white] (0,0) --(315:2.5)--(360:2) ;
     \fill(0,0) node[sing] {};
\end{tikzpicture}
    \caption{Canonical scattering diagram}
    \label{fig:affinescatter}
\end{figure}

Now we would like to define the canonical scattering diagrams from $(B, \Sigma)$. 
Rather than obtaining the diagrams by the algorithmic process with some initial data in \cite{GS} \cite{KS}, the canonical scattering diagrams are defined via some Gromov-Witten type invariants.
We will discuss the two types of diagrams are the `same' later in the discussion about how to go from canonical scattering diagrams to cluster scattering diagrams. 
A \emph{wall} \cite{ghks_cubic} in $B$ is a pair $(\fd, f_{\fd})$ where $\fd \subset \sigma_{i,i+1}$, for some $i$, is a ray generated by $a v_i+ bv_{i+1} \neq 0$, $a, b \in \Z$ relatively prime, and  $f_{\fd} = 1 + \sum_{k \geq 1} c_k X_i^{-ak}X_{i+1}^{-bk} \in \C[[X_i^{-a}X_{i+1}^{-b}]]$ with some  finiteness properties, and where $c_k $ corresponds to the curve counting invariants.
Note that the description of the wall functions $f_{\fd}$ indicates that all the wall are outgoing in the sense stated in the last section.
Then the scattering diagrams are again the collections of walls. 
For example, the canonical scattering diagram associated to Figure \ref{fig:affine} is shown in Figure \ref{fig:affinescatter}.

Let $B(\Z)$ be the set of points of $B_0$ with integral coordinates in an integral affine chart and $\{ 0\}$. 
Theta functions $\vartheta_q$, $q \in B(\Z)$, can similarly be defined on the scattering diagrams. 
The set of theta functions again generates an algebra structure \cite{GHKlog} in terms of broken lines.
In the finite case, we can simply consider $A = \oplus_{q \in B(\Z)} \vartheta_q$. 

Analogous to the setting in the cluster scattering diagrams, we can use the Rees construction to compactify the mirrors \cite{gross2019intrinsic}.
Positive polytopes with respect to the affine structures can be similarly defined to give a graded algebra. 
Using \cite{mandel2016tropical} or the argument in \cite{cpt}, the polytopes are broken line convex. In this case, since all the walls are outgoing, the positive polytopes are simply convex with respect to the affine structures.

\subsubsection*{Relation to the cluster scattering diagrams}

In the case of $Y$ a non-singular toric surface and $D= \partial Y$ the toric boundary of $D$, the affine structure on $B $ extends across the origin.
This identifies $(B, \Sigma)$ with $M_{\R}, \Sigma_Y$, where $\Sigma_Y$ is a fan for $Y$. 

Now given a Looijenga pair. 
Assume there is a toric model $p: (Y, D) \rightarrow ( \bar{Y} , \bar{D})$ which blows up distinct points $x_{ij}$ on $D_i$.
A \emph{toric model} of $(Y, D)$ is a birational morphism $(Y, D) \rightarrow (\bar{Y}, \bar{D})$ to a smooth toric surface $\bar{Y}$ with its toric boundary $\bar{D}$ such that $D \rightarrow \bar{D}$ is an isomorphism. 
Consider the tropicalisation $(\bar{B}, \bar{\Sigma})$ of $(\bar{Y}, \bar{D})$. 
Thus $\bar{B} \cong M_{\R} = \R^2$ and $\bar{\Sigma}$ is the fan for $\bar{Y}$. 
Then there is a canonical piecewise linear map
\[
\nu : B \rightarrow \bar{B}
\]
which restricts to an integral affine isomorphism on the maximal cones in $\sigma$ and $\bar{\Sigma}$. 
One can then define the scattering diagram $\fD$ as outlined in Section \ref{sec:ClusterVar} or as in \cite[Definition 3.21]{GHKlog} for the more general setting. 
This step can be seen as `pushing the singularities to infinity' or `moving worms' \cite{KS}. 
By definition, the singularity of the affine structure is at $\{0\}$ as indicated in Figure \ref{fig:affine}. 
Then the singularity can be imagined to be pushed to the infinity of the two incoming walls. 
The map $\nu$ can be extended to act on the canonical scattering diagram $\fD^{\text{can}}$. 
It is shown that \cite{GHKlog} $\bar{\fD} = \nu (\fD^{\text{can}})$. 

We have discussed in Section \ref{sec:ClusterVar} that every cluster variety can be described as blow ups of a toric variety, which give the toric models for the cluster variety. 
Thus the cluster scattering diagrams can be seen as the diagrams arising from the canonical scattering diagrams by the map $\nu$ (as pushing singularities to infinity).

\subsubsection*{Mutation of polytopes according to the affine structures}

\begin{figure}
\centering
    \begin{tikzpicture}
    \fill(0,0) node[sing] {};
    \fill(0:3) node[sing] {};
    \draw[->] (0) -- ++(180:2) node[left]{$1+z^{e_1}$}; 
    \draw[->] (0)--++(90:2) node[above]{$1+z^{-e_2}$}; 
    \draw[->] (0)--++(225:2) node[left]{$1+z^{e_1+e_2}$};
    \draw (0)--++(0:2) ;
    \draw[<-] (0:1)--++(0:1);
    \draw[->] (0)--++(270:2) node[below]{$1+z^{e_2}$}; 
     \draw[fill=gray, opacity=0.2, draw=white] (0,0) --(135:2.5)--(90:2) ;
\end{tikzpicture} 
     \caption{Scattering diagrams with monodromy.} \label{fig:monommoo}
\end{figure}
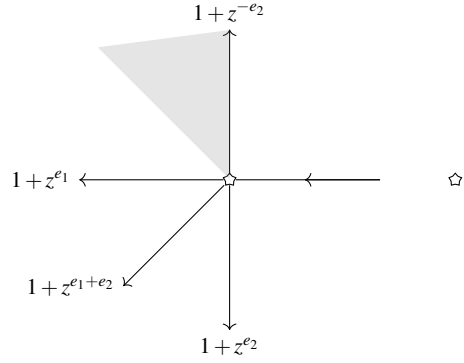

\begin{figure}[h!]   
  \begin{center} 
 \centerline{\includegraphics[scale=3]{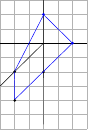}}
\caption{Positive polytope with respect with the underlying affine structure} \label{fig:polytopewmono}
\end{center}
\end{figure}

One can imagine or with symplectic motivation as in Section \ref{sec:SYZ}, the `pushing singularities to infinities' procedure is more general than just having singularities at the origin.
For example, we can consider Figure \ref{fig:monommoo} which we only push one of the singularities to infinity, resulting in a scattering diagram with one incoming wall. 
The monodromy in Figure \ref{fig:monommoo} is 
$\begin{pmatrix}
	1 & 0 \\ 1 & 1
\end{pmatrix}$, or $(1,0) \mapsto (1,1)$, $(0,1) \mapsto (0,1)$. 
The polytope in Figure \ref{fig:A2polytope} with respect to this affine structure would then be of the form in Figure \ref{fig:polytopewmono}. 

\begin{figure}[h!]   
  \begin{center} 
 \centerline{\includegraphics[scale= 2.2]{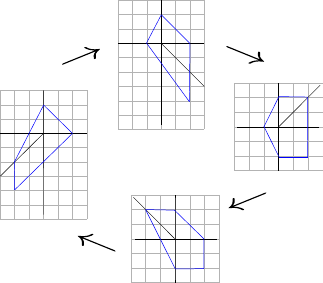}}
\caption{Mutation of polytopes with monodromy} \label{fig:seq}
\end{center} 
\end{figure}

We can apply the sequence of mutations in Figure \ref{fig:XA2mutate} to the polytope in Figure \ref{fig:polytopewmono} and then obtain a new sequence of  mutation polytopes (Figure \ref{fig:seq}). Putting the polytope in Figure \ref{fig:A2polytope} into the sequence  (Figure \ref{fig:seq}) will get us the sequence Figure \ref{fig:dp5_1} which is motivated from the symplectic perspective. 

\begin{figure}
\centering
    \begin{tikzpicture}
    \draw[->] (1,0) -- (2,0) node[right]{$1+z^{-e_1}$}; 
    \draw[->] (0,1)-- (0,-2) node[below]{$1+z^{e_2}$}; 
    \draw[->] (0,0) -- (-2,-2) node[left]{$1+z^{e_1+e_2}$};
    \draw[->] (1,0) -- (-2,0) node[left]{$1+z^{e_1}$};
    \draw[->] (0,1) -- (0,2) node[above]{$1+z^{-e_2}$}; 
     \draw[fill=gray, opacity=0.2, draw=white] (1,0) --(2,1)--(2,0) ;
         \fill(1,0) node[sing] {};
    \node[right] at (2,1) {$\begin{pmatrix}
	1 & -1 \\ 0 & 1 \end{pmatrix}$};
	\draw[fill=gray, opacity=0.2, draw=white] (0,1) --(-1,2)--(0,2) ;
         \fill(1,0) node[sing] {};
	\node[left] at (-1,2) {$\begin{pmatrix}
	1 & 0 \\ 1 & 1 \end{pmatrix}$};
    \fill(0,1) node[sing] {};
\end{tikzpicture} 
     \caption{Scattering diagram with monodromy on the walls.} \label{fig:monoD}
\end{figure}
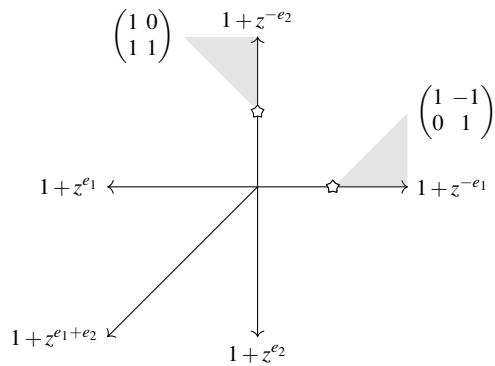

The singularities can also be located on the walls instead of just at infinity or the origin. 
For example, one can obtain the canonical scattering diagram shown in Figure \ref{fig:monoD}. 
Similar calculation shown in \cite{dp5} indicates that the scattering diagram is consistent.

\begin{figure}[h!]   
  \begin{center} 
 \centerline{\includegraphics[scale= .6]{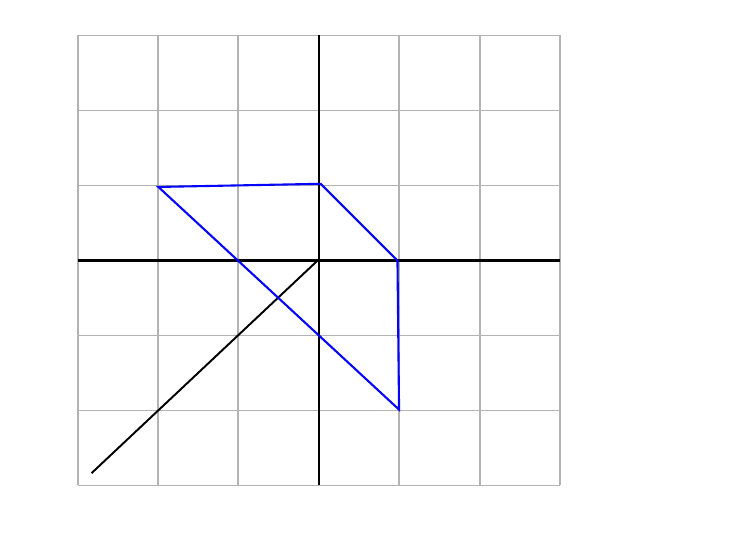}}
\caption{Polytope lives in the affine structure indicated in Figure \ref{fig:monoD}.} \label{fig:monopoly}
\end{center} 
\end{figure}

Note that the portions of the walls which go from the singularities to infinity are all outgoing. 
Thus using the idea in \cite[Remark 6.2]{cpt}, we can consider convex sets in this affine structure. 
For example, one can construct the polytope as in Figure \ref{fig:monopoly}. 

\begin{figure}[h!]   
  \begin{center} 
 \centerline{\includegraphics[scale= 2.5]{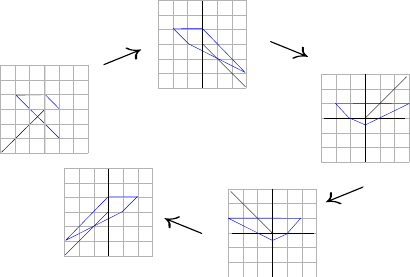}}
\caption{Mutation for the polytope in Figure \ref{fig:monopoly}} \label{fig:monopoly_mutate}
\end{center} 
\end{figure}

Applying the mutation sequence as in Figure \ref{fig:XA2mutate},
we obtain the sequence of polytopes described in Figure \ref{fig:monopoly_mutate}.
Interestingly, this is the same sequence as in Figure \ref{fig:dp8_1} which is motivated from the symplectic perspective. 

\subsection{Type $B_2$ and $G_2$}

The other types are similar and thus we only roughly go over the mutations of type $B_2$ and $G_2$. 
For type $B_2$, we will take the same skew-symmetric form with $d_1 =1$, and $d_2=2$ as our fixed data. We can again take the initial seed as $\seed = \{ (1,0), (0,1)\}$.
Then we will obtain the $\cA$ and $\cX$ scattering diagrams as in Figure \ref{fig:rank2} and \ref{fig:Xrank2}. 
If we mutate at index $1$ first, we can get the mutation of scattering diagrams very similar to the type $A_2$ case.

\begin{figure}[H]   
  \begin{center} 
 \centerline{\includegraphics[scale= 0.4]{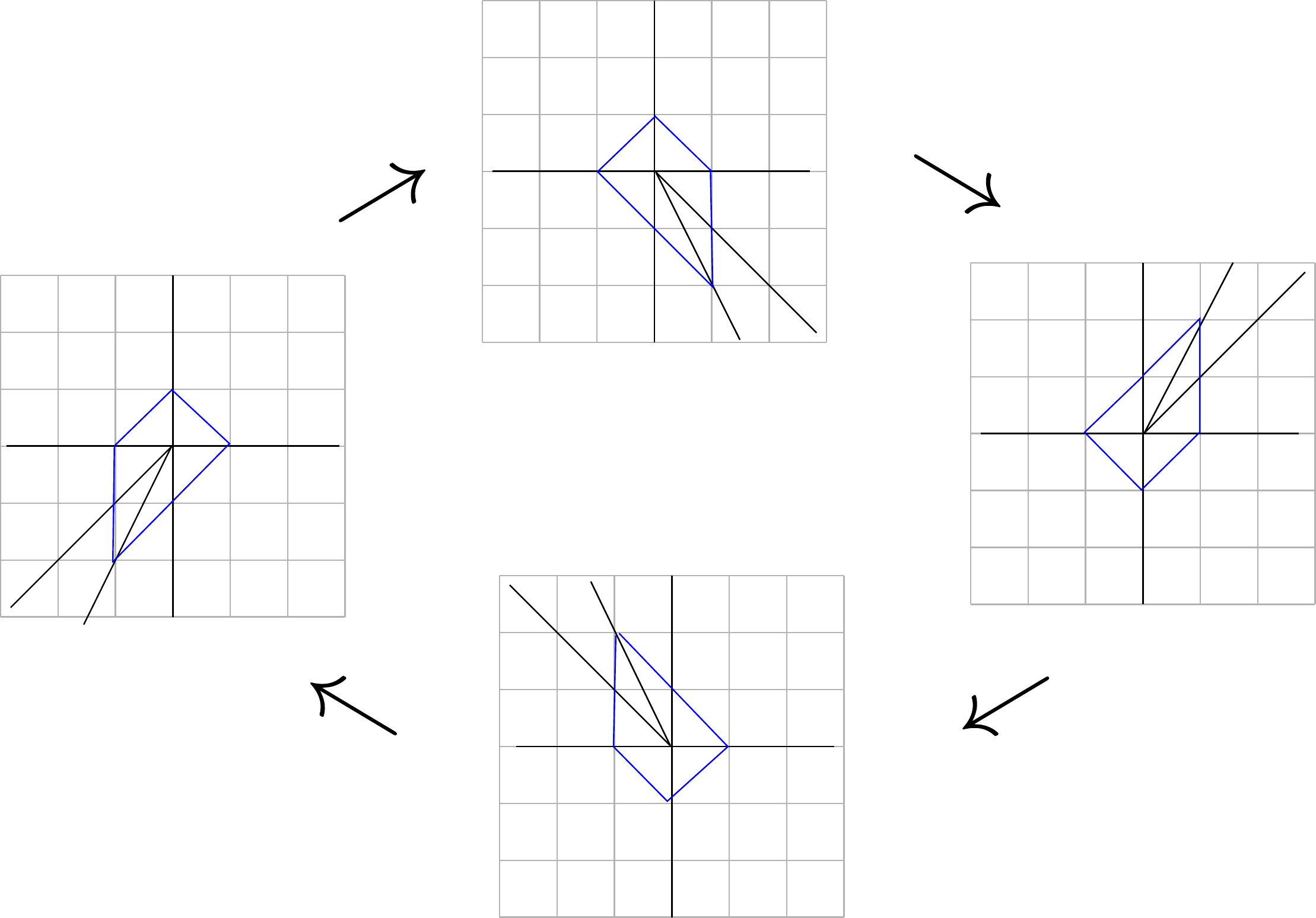}}
\caption{Mutation of polytope for type $B_2$} 
\label{fig:XB2mut_polytope}
\end{center} 
\end{figure}

For type $B_2$, we can again take the primitive generators of the walls and then consider the polytope as the convex hull of those vertices. 
By using \cite{cpt}, this polytope is a positive polytope. 
The multiplication of the theta functions tells us \cite{bat_finitetype} that the corresponding space is the del Pezzo surfaces of degree 6. 
The mutation sequence of the polytopes is described in Figure \ref{fig:XB2mut_polytope}.

One may want to repeat the same trick on the type $G_2$. 
The sad fact is that the if we are taking the convex hull of the primitive generators of the walls, the resulting polytope would no longer be positive. 
This is because the polytope is no longer broken line convex as indicated in \cite{cpt}. 
Since we only care about the incoming walls for broken line convexity \cite{cpt}, one can see that the top left polytope in Figure \ref{fig:dp6_1} in the next section is broken line convex. 
The mutation sequence of the $G_2$ scattering diagrams are similar to those for type $A_2$ and $B_2$. 
Without duplicating, one can see that the mutation indicate in Figure \ref{fig:dp6_1} is in fact the cluster mutation of scattering diagrams. 
Thus the mutation of polytopes follows correspondingly. This again tells us that the mutation sequences of the polytope with the algebro-geometric and symplectic viewpoints coincide.

\section{Mutations in symplectic geometry} \label{sec:symp}

We begin this section giving a perspective on understanding cluster varieties, as well as
scattering diagrams, as a way of building mirrors under SYZ \cite{SYZ96} $T$-duality. 
In particular, we explain how scattering diagrams can be related to almost toric fibrations. Later, we explain how one can see compactifications of 
$2$-dimensional cluster varieties into del Pezzo surfaces from an almost toric fibration perspective.

\subsection{Cluster varieties and Mirror Symmetry}\label{sec:SYZ}

In this section we will sketch how to relate a scattering diagram data (described in Section \ref{sec:ag}), for constructing a log-CY variety $X$, 
with the base of an almost-toric fibration (ATF), describing a SYZ \cite{SYZ96} singular Lagrangian fibration of $(X,\omega)$, 
with respect to a K\"ahler form $\omega$ in $X$. 

\subsubsection{Almost Toric Fibrations}

Informally speaking an \emph{ almost toric fibrations} (ATF) in a symplectic 4-manifold $X$ is a smooth map to a two dimensional base $B$, whose regular fibres are Lagrangian tori, whose allowed singular fibres are of three kinds: 

\begin{itemize}
    \item point (toric - rank 0 elliptic) -- locally equivalent to the moment map at the origin in $\C^2$ with the standard toric action, $(e^{i\theta_1},e^{i\theta_2})\cdot(x,y) = (e^{i\theta_1}x,e^{i\theta_2}y)$;
    \item circle (toric - rank 1 elliptic) -- locally equivalent to $S^1\times \{0\} \subset \C^* \times \C$ with the standard toric action;
    \item nodal  (a pinched torus) -- with some local model described for the singular point. [See \cite{Sy03,SyLe10} for precise definition, and see Section~\ref{subsec:LocModel} for a local model of the nodal fibre.]
\end{itemize}
   The toric singularities appear on the boundary of the base, 
while the nodal singularities project into the interior. For a precise definition of ATFs see \cite{Sy03}.

Away from the singular fibres, by the Arnold-Liouville theorem \cite{Ar_book}, $X$ admits locally action angle coordinates $
(p_1,p_2,\theta_1,\theta_2)$ and the fibration is locally equivalent to $
(p_1,p_2,\theta_1,\theta_2) \mapsto (p_1,p_2)$, in other words, away from singular fibres $X$ equivalent to $T^*B/ \Lambda^*$, 
for some lattice $\Lambda^*$. Hence, $B$ carries a natural dual lattice $\Lambda \subset TB$. The lattice has monodromy
as we go around the nodal fibre, which is a shear in the direction dual to the collapsing cycle of that nodal fibre.
Locally, the coordinates $(p_1,p_2)$, can be thought as the \emph{flux} $\ff \in H^1(T^2,\R)$ relative to the Lagrangian fibre associated with $(0,0)$. The flux $\ff(\gamma)$ measures the symplectic area of a cylinder swept by 
a cycle $\gamma \in  H_1(T^2,\Z)$ as we move in a path of Lagrangian fibres connecting  $(0,0)$ to $(p_1,p_2)$.
[See, for instance, \cite{ShToVi18} for a more complete understanding of flux in ATFs.]

So, in practice, we visualise the base minus a set of cuts (one for each nodal fibre) affinely embedded into $\R^2$ endowed with the standard affine structure. We call them almost-toric base diagrams (ATBDs)
representing the ATF. The same ATF can be represented by different ATBDs, by changing the set of cuts. 

\begin{figure}[h!]
\sidecaption
\centerline{\includegraphics[scale= 1]{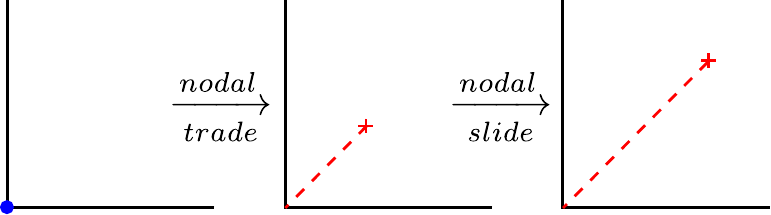}}
\caption{Nodal trade and nodal slide operations in ATFs.}
\label{fig:ATF_ops} 
\end{figure}

Figure \ref{fig:ATF_ops} shows the base diagram of 3 different ATFs in $\C^2$; the right-most 
diagrams on Figure \ref{fig:SYZ2} are different diagrams representing the same ATF in $\C^2 \setminus \{xy=1\}$, 
related by a change of cut; Figures \ref{fig:dp5_2_1}--\ref{fig:dp4_1} contain examples of ATBDs in closed 4 manifolds.
In these diagrams, the crosses represent the nodal fibres, 
the dashed lines the cuts, the edges the rank 1 and the dots rank 0 toric singularities. 

\begin{remark}
We expect the above mentioned ATFs to be realisable as a special Lagrangian fibration in the complement of a 
complex divisor projecting to the boundary of the ATF, with respect to a holomorphic volume form with poles on these divisor. 
This is true for the fibration presented in Section \ref{subsec:LocModel}, but we will avoid talking about the "special" condition.
\end{remark}

There are two ways of modifying ATFs within the same symplectic manifold $X$, known as nodal trade and nodal slide \cite{Sy03}. 
The diagrams in Figure \ref{fig:ATF_ops} illustrate the change of the ATBDs after a nodal trade and a nodal slide. 
In del Pezzo surfaces, the monotone symplectic 4-manifolds, we defined mutation of an ATBD, the process of sliding one nodal
fiber through the monotone fibre, and then redrawing the diagram by changing the direction of the cut used to slide \cite{Vi14,Vi16a}. 
In the end, the ATBD mutates by slicing it in the direction of the cut and applying the inverse of the corresponding monodromy, 
which is a shear in the primitive direction associated to the cut. This transformation is the same polytope mutation as in \cite{ACGK12,AkKa16}, 
and completely analogous to the mutation of seeds, and scattering diagram we will discuss later.
We can extend this notion of symplectic mutation to exact almost toric manifold, 
for instance, the complement of an anti-canonical divisor in a del Pezzo. In this case, the mutation corresponds to
sliding a nodal fibre through the exact torus and then transferring the associated cut to the opposite side. 

\subsubsection{Local model for nodal fibre and wall-crossing} \label{subsec:LocModel}

We briefly recall ATF presented in \cite[Section~5]{Au07}, \cite[Section~3.1.1]{Au09}. This ATF appeared before in 
\cite{ElPo93} and also in \cite[Example~1.2]{Gro97}, where it was shown to be a special Lagrangian fibration [with respect
to certain holomorphic volume form].
We consider $X^\vee = \C^2 \setminus \{xy = 1\}$, with $\omega^\vee = \frac{i}{2} \left( dx\wedge d\bar{x} + dy\wedge d\bar{y} \right)$ the standard symplectic form. 
Using $f: X^\vee \to \C \setminus \{1\} $, $f(x,y)=xy$, Auroux builds an ATF by parallel transport of orbits of the $S^1$ action $e^{i\theta}\cdot(x,y) = (e^{i\theta}x,e^{-i\theta}y)$, over circles in the base of $f$ centred at $1$. One then gets Lagrangian torus fibres, parametrised by $(r,\lambda) \in \R_{>0}\times \R$, as:
\[ T_{r,\lambda} = \{(x,y) \in \C^2; r = |xy - 1|, \lambda = |x|^2 - |y|^2\}. \]

Note that there is a nodal fibre $T_{1,0}$, that contains $(0,0)$, the fixed point of the $S^1$ action.

\begin{figure}[h!]   
  \sidecaption
  \centerline{\includegraphics[scale= 0.7]{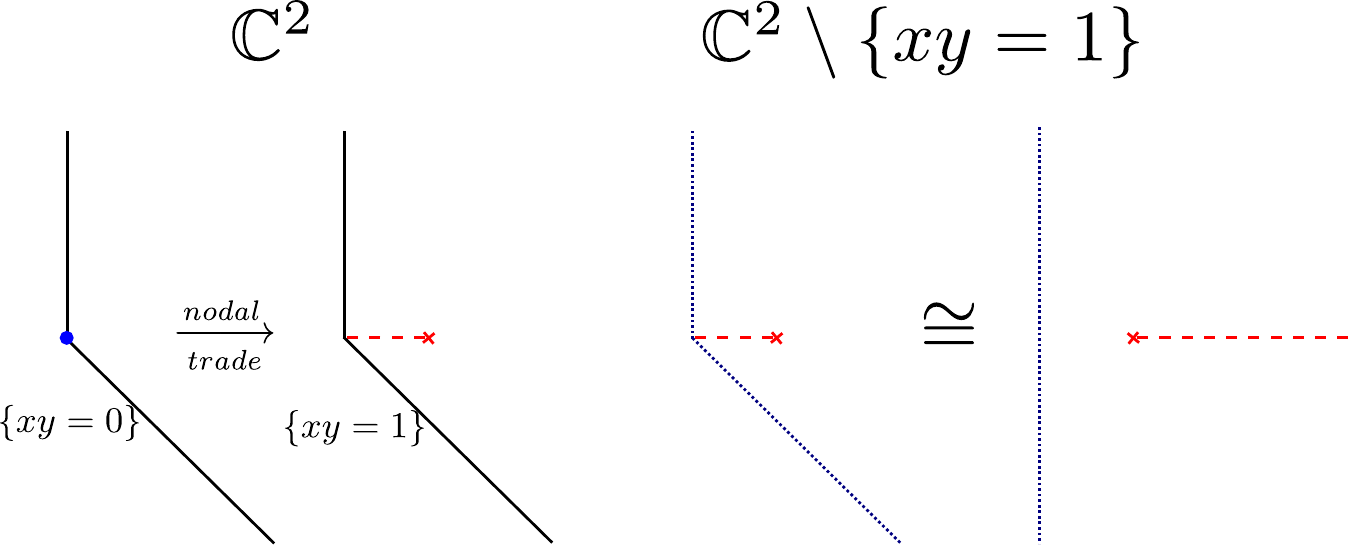}}
\caption{Nodal trade and an ATF for the complement of a conic. The left diagrams are a shear by $(0,-1)$ of the diagrams in 
Figure \ref{fig:ATF_ops}. The rightmost diagrams, represent the same ATF, and differ by changing the direction of the cut.}
\label{fig:SYZ2} 
\end{figure}

This almost toric fibration can be represented by applying a nodal trade to the standard toric fibration of $\C^2$, replacing the
boundary divisor $\{xy = 0\}$ with the smooth divisor $\{xy = 1\}$, and then deleting this divisor living over the 
boundary of the base, as illustrated by Figure \ref{fig:SYZ2}. Indeed, replacing the role of $1$ by $0$ in the above fibration, i.e., 
considering parallel transport over circles concentric at $0$ (considering $r = |xy|$) one obtain precisely the standard toric
fibration of $\C^2$. So, considering analogous fibrations by changing $0$ to $1$ in the definition of $r$ constitutes a nodal trade, and moreover, varying the value of $c \in \R_{>0}$ in the definition of $r = |xy - c|$ provides different fibrations related by nodal slides. 

\subsubsection*{Wall-crossing}

Dualising this torus fibration, one gets the mirror variety $X_\Lambda$ of $X^\vee$, over the
Novikov field $\Lambda = \{ \sum_{i=0}^n a_iT^{\sigma_i} ; a_i \in \C, \sigma_i \in \R, \lim_i \sigma_i = \infty \}$, 
which is the moduli of almost-toric fibres (special Lagrangians), endowed with unitary $\Lambda^*$-local systems.
We will be able to relate the valuation $\val (u) = \min \{ \sigma_i; a_i \ne 0\}$ of an element $u \in \Lambda$ with the
above mentioned flux, 
whenever $\val (u)$ measures the symplectic area of a disk with boundary in a varying family of the Lagrangian torus fibres. 
[Notation: $\Lambda_0 = \{ u \in \Lambda; \val (u) \ge 0 \}$, $\Lambda_+ = \{ u \in \Lambda; \val (u) > 0 \}$,
$\Lambda^* = \{a_0 + \Lambda_+; a_0 \in \C^*\}$.]
We will later consider the mirror of $X^\vee$ as
$X = X_\C$ over $\C$, by replacing $T$ with $e^{-1}$. 

So we replace the Lagrangian fibre $T^2$, by the dual $\Lambda$-torus of unitary local systems 
$\hom(\pi_1(T^2); \Lambda^*) \cong ( \Lambda^*)^2$. Locally identifying each relative class, 
$\beta \in \pi_2(\C, T_{r,\lambda})$,
Auroux defined a function $z_\beta: X_\Lambda \to \Lambda$, for a local system $\nabla$ in $T_{r,\lambda}$, 
$z_\beta(\nabla) = T^{\omega^\vee(\beta)}\nabla \cdot \del \beta$. 
Choosing a basis $\{\alpha,\beta\}$ of $\pi_2(\C, T_{r,\lambda})$, one gets that $w  := z_\alpha, \ u := z_\beta$, 
define local coordinates of $X_\Lambda$. After that, the idea is to define a superpotential function 
$W: X_\Lambda \to \Lambda_+$, which is locally defined as $W(u,w)$, and whose monomials encode the relative Gromov-Witten 
count of Maslov index 2 holomorphic disks in $\C$ with boundary on the torus fibre endowed with the respective local system 
determined by $(u,w)$. [The pair $(X_\Lambda, W)$ is called the Landau-Ginzburg model that is
mirror dual to $\C$ with respect to the divisor $D = \{xy - 1\}$. We refer the reader to \cite{Au07,Au09} for details on 
mirror symmetry in the complement of divisors.] 

The issue is that, in the naive definition of the mirror, the superpotential $W$ is discontinuous. This is due to the presence of fibres $T_{1,\lambda}$, 
$\lambda \ne 0$, which bounds Maslov index $0$ holomorphic disks.  Let's denote the relative class represented by this Maslov 0 disks
by $\alpha$ for $\lambda <0$, and $-\alpha$ for $\lambda>0$. In \cite[Section~3.1.1]{Au09}, it is shown that
for $r < 1$, the fibres $T_{r,\lambda}$ (called Chekanov type) bound one holomorphic disk, in a class we name $\beta$. So 
$W(u,w) = u$, for these fibres. The fibres $T_{r,\lambda}$, for $r > 1$, (called Clifford type) bound 2 holomorphic disks in
relative classes $\beta_1$, $\beta_2$, and hence the superpotential is of the form $W(z_1,z_2) = z_1 + z_2$, where $z_i =z_{\beta_i}$. 

We see in  \cite[Section~3.1.1]{Au09} that as $r$ approaches $1$, from $r > 1$, we get $\alpha = \beta_1 - \beta_2$
(hence $w = z_1z_2^{-1}$).
Moreover, if we cross the wall at $\lambda < 0$, the class $\beta$ is naturally identified with $\beta_2$, and 
 if we cross the wall at $\lambda > 0$, the class $\beta$ is naturally identified with $\beta_1 = \beta_2 + \alpha$
 [which is not so surprising, as the monodromy around the nodal fibre in the ATF would fix $\del \alpha$ and maps 
 $\del \beta \to \del \beta + \del \alpha$]. So the superpotential $W$ should be corrected by the term $(1 + w^{\pm1})$,
 representing the fact that the holomorphic disk on class $\beta$, would not only survive past the wall, but the superpotential
 would also acquire a holomorphic disk in class $\beta \pm \alpha$, coming from the gluing of the Maslov 2 holomorphic disk on 
 class $\beta$ with the Maslov 0 holomorphic disk on class $\pm \alpha$.
Then, instead of $u$ becoming $z_2$ as we cross over $\lambda > 0$,  we should correct it to become $u = z_2(1+w) = z_2 + z_1$, and instead of $u$ becoming $z_1$ as we cross over $\lambda < 0$, we should correct it to become $u = z_1(1+w^{-1}) = z_1 + z_2$, 
and, thus, ensuring the continuity of $W$. 

By naming $v = z_2^{-1}$, so $z_1 = v^{-1}w$, we get the corrected $u = v^{-1}(1+w) = v^{-1}w(1+w^{-1})$. We see that the corrected (and completed) mirror $X_\Lambda$, is given by 

\[ X_\Lambda = \{(u,v,w) \in \Lambda^2\times (\Lambda\setminus \{0\}); uv = 1 + w \}.\] 

\begin{figure}[h!]   
  \sidecaption
  \centerline{\includegraphics[scale= 0.5]{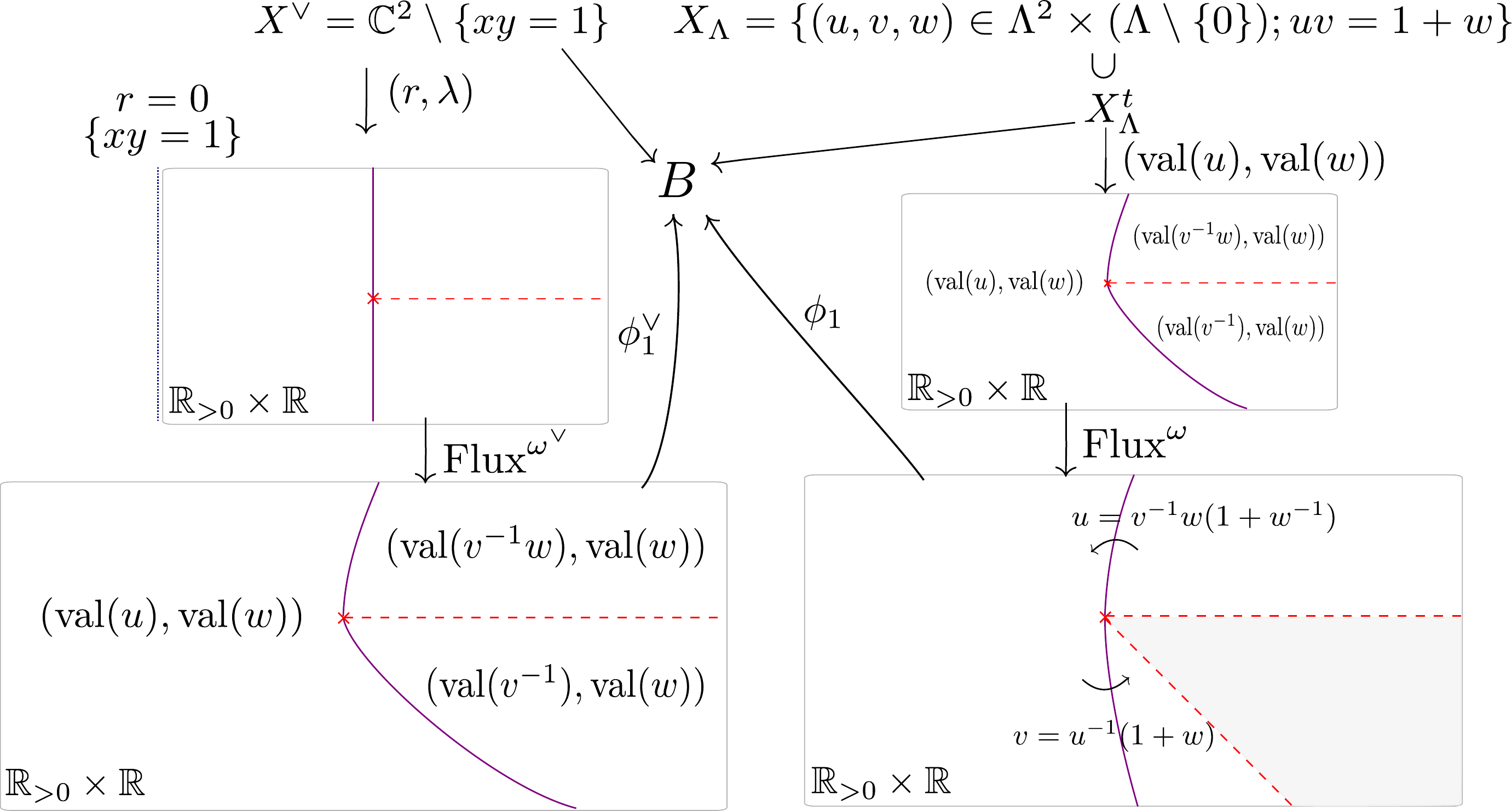}}
\caption{SYZ fibrations for the complement of a conic, which is self-mirror, when considering $X_\C$.}
\label{fig:SYZ1} 
\end{figure}

Figure \ref{fig:SYZ1} below describes the \emph{mirror  SYZ fibrations} on $X^\vee$ and $X_\Lambda$. We list 
several remarks about the diagrams in Figure \ref{fig:SYZ1} and the mirror $X_\Lambda$.

\begin{remark}
We see that $X_\Lambda$ is given by gluing two torus charts $(u,w) \in (\Lambda\setminus \{0\})^2$, and 
$(v,w) \in (\Lambda\setminus \{0\})^2$, by a rational map defined in the complement of $\{w = -1\}$. So, for instance, 
the case $u=0$ would be realised as $(v,-1)$ in the $(v,w)$-chart.  
\end{remark}

\begin{remark}
Considering the symplectic form $\omega^\vee$ in $X^\vee$, the Lagrangian torus fibration would be viewed 
in a truncated part of the $(u,w)$-chart, with $0 < \val(u) \le a(\val(w))$ or in the $(v,w)$-chart, 
with $a(\val(w)) \le \val(v^{-1}) < \infty$ for $\val(w) < 0$, for instance. These bounds on the valuation would give us $X_\Lambda^t$, a truncated version of the mirror $X_\Lambda$. But in symplectic geometry, we can add to $(X^\vee, \omega^\vee)$ 
a contact boundary $\del X^\vee$, and it is most natural to consider a completion procedure called the symplectization of $X^\vee$ with respect to this 
boundary. This endows $X^\vee$ with a different symplectic form $\omega^\vee_S$. It is equivalent to consider an infinite
inflation of $(\C^2,\omega^\vee)$ with respect to the divisor $D = \{xy = 1\}$. In this limit we would have
$\val(u) \to \infty$, and we would get the completed mirror $X_\Lambda$.
\end{remark}

\begin{remark} The expectation regarding the correspondence between the count of Maslov index 2 disks with boundary on a SYZ fibre and
its tropical counterpart was proven in \cite{lin2020enumerative} for a SYZ fibration on the complement of an smooth anti-canonical divisor 
in a del Pezzo surface.
More precisely, given a del Pezzo surface $Y$ and a smooth anti-canonical divisor $D$, there exists a special Lagrangian fibration on $Y\setminus D$ with respect to the complete Ricci-flat Tian-Yau metric \cite{CJL}. To understand the Landau-Ginzburg superpotential of $Y$, there exists a sequence of K\"ahler forms $\omega_i$ on $Y$ converging to the Tian-Yau metric pointwisely with $\int_Y \omega_i^2\rightarrow \infty$ \cite[Lemma 2.4]{lin2020enumerative}. Thus, these superpotentials of the special Lagrangian fibres can be defined with respect to $\omega_i$, $i \gg 0$ and the superpotentials coincide with the 
tropical counterpart \cite[Theorem 5.19]{lin2020enumerative}. This gives a geometric explanation of the renormalization procedure of taking valuation going to infinity.
\end{remark}

\begin{remark}
In the $(r,\lambda)$ projection of Figure \ref{fig:SYZ1}, the singular fibre in position $(1,0)$ is depicted 
by an $\times$, and the wall of fibres with $r=1$ that bound Maslov index 0 disks are represented by a line.
This $(r,\lambda)$ coordinate does not respect the natural affine structure on the complement of the
singular fibre of $B$. We instead consider $\mathrm{Flux}^{\omega^{\vee}}$, the flux with respect to a limiting fibre lying over $(0,0)$, in the next diagram. 
This map is then continuous, but not differentiable over the dashed ray $r\ge 1$, $\lambda = 0$, which we call the cut. Moreover, 
this composition represents the map to the base diagram depicted in the rightmost picture of Figure \ref{fig:SYZ2}.
The map $\phi_1^\vee$ is then an affine isomorphism to $B$ minus the cut. The affine structure of $B$ minus the node, 
is described by the gluing of the chart $\phi_1^\vee$ and a chart $\phi_2^\vee$, going from the third diagram of Figure
\ref{fig:SYZ2}, corresponding to taking the cut associated to $0 < r \le 1$, $\lambda = 0$. The symplectic manifold
$X^\vee$ can be thought then as a local model for gluing in the nodal fibre to the manifold constructed from gluing
the Lagrangian torus fibres associated to $\phi_1^\vee$ and $\phi_2^\vee$. This is essentially the same model
as the description of $X^\vee$ as a self-plumbing of $T^*S^2$ given in \cite[Section~4.2]{Sy03}.

\end{remark}

\begin{remark}
The affine structure on $B$ for the dual mirror fibration $X_\Lambda \to B$ is endowed with the dual affine structure in the complement of the node. We call the map that adjusts this affine structure in $\R^2$, $\mathrm{Flux}^{\omega}$. In the case we
take the SYZ mirror $X_\C$ over $\C$ (by replacing $T$ by $e^{-1}$), it endows a symplectic form $\omega$ as described in 
\cite[Proposition~2.3]{Au09}. In this case, $\mathrm{Flux}^{\omega}$ becomes the actual flux with respect to this symplectic form. 
\end{remark}

\begin{remark}
The monodromy around the singular fibre of $X^\vee \to B$, represented by the bottom left diagram of Figure 
\ref{fig:SYZ1}, is given by $M^{\pm 1}$, for $M =\begin{bmatrix} 1 & -1 \\ 0 & 1 \end{bmatrix}$, fixing the cut $(1,0)$.
Then, the monodromy around the node for the rightmost diagram representing $X_\Lambda \to B$ is given by 
$(M^T)^{\mp 1}$, with $(M^T)^{- 1} =\begin{bmatrix} 1 & 0 \\ 1 & 1 \end{bmatrix}$. We see that this fixes the 
coordinate $w$, associated to $(0,1)$, which we then name $\vartheta_{(0,1)}$, and it sends the coordinate
$v^{-1} = \vartheta_{(1,0)}$ to $v^{-1}w = \vartheta_{(1,1)}$. 

\end{remark}

\begin{figure}[h!]   
  \sidecaption
  \centerline{\includegraphics[scale= 0.5]{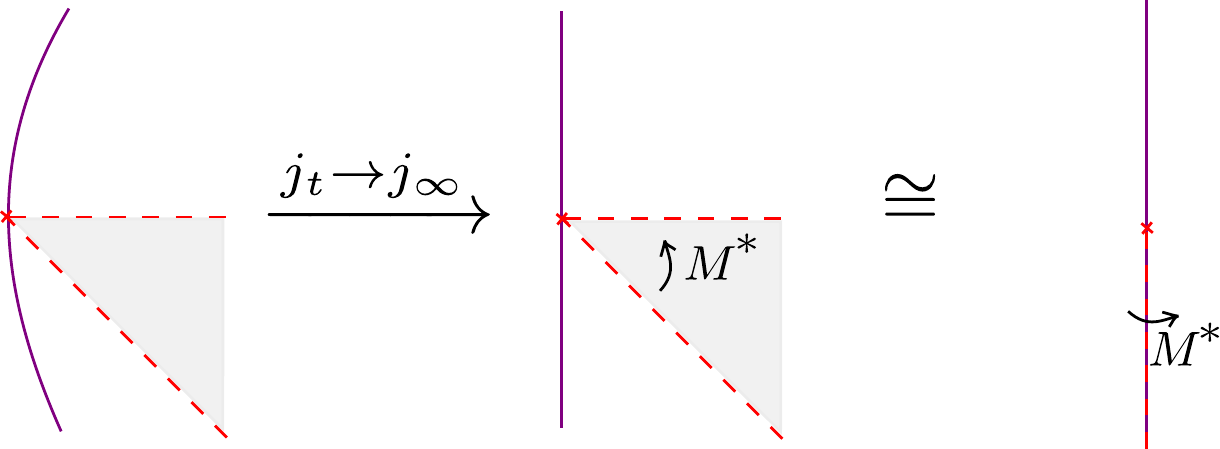}}
\caption{In a complex structure limit, the wall becomes straight. We can move the cut to the invariant direction of the monodromy,
that is the same direction as the limit straight wall.}
\label{fig:SYZ3} 
\end{figure}

\begin{remark}
As described by Mikhalkin \cite{Mi04b}, we can deform the complex structure on $X^\vee$ to a limit where holomorphic curves
would converge to tropical curve on the base with respect to the so-called complex affine structure, which is dual to the symplectic
affine structure. So, (relative) 
Gromov-Witten invariants of $X^\vee$ are expected to be described by tropical curves ($\vartheta$ functions) in the base $B$, 
with the affine structure describing $X_\Lambda$ (or $X_\C$). In particular, the wall becomes straight in this limit, 
as illustrated by Figure \ref{fig:SYZ3}. 
\end{remark}

As we mentioned before, we will now replace the formal variable $T$ by $e^{-1}$ in our construction, and consider the
mirror as the moduli of Lagrangian fibres endowed with $U(1)$-local systems. So, after completion the mirror becomes 

\[ X = X_\C = \{(u,v,w) \in \C^2\times \C^* ; uv = 1 + w \} \] 
endowed with a completed symplectic form $\omega$ as described in \cite[Proposition~2.3]{Au09}. Its SYZ dual ATF (dual to the one
on $X^\vee$) is then described by any of the diagrams in Figure \ref{fig:SYZ3}. 

\begin{remark}
As we take the completion, the $\val$ defined in $X_\Lambda$ approaches $-\log|.|$ defined in $X_\C$.   
\end{remark}

We see now that the rightmost diagram in Figure \ref{fig:SYZ3}, can describe an ATF, and once decorated with 
wall crossing functions [$(1 + w^{\pm1})$ accordingly] along the wall, it can algebraically determine the space
$X_\C$. The variety $X_\C$ is then built out of two $(\C^*)^2$ charts, with coordinates $(u,w)$ and $(v,w)$, 
glued together in a cluster like transition birational map $uv = 1 + w$, defined in the complement of $\{w = -1\}$. 
A unique wall, decorated with such wall crossing function, describing $X_\C$  is the simplest version of a scattering diagram \cite{GS,GHKlog} (see Section \ref{sec:ag} for more details).


\subsubsection{The $A_2$ Cluster Variety: ATF and Scattering Diagram} \label{sec:A2}

As described in \cite[Example~8.40]{ghkk} (taking the parameters $X_1, X_2$ to be $1$), 
we describe the affine ($\vartheta_0 = 1$) $A_2$ Cluster variety by the ring in five variables $\vartheta_i$, 
$i = 1,\dots, 5$ satisfying the relations:

\[ \vartheta_1\vartheta_3 = 1 + \vartheta_2 \]
\[ \vartheta_2\vartheta_4 = 1 + \vartheta_3 \] 
\[\vartheta_3\vartheta_5 = 1 + \vartheta_4 \]
\[\vartheta_4\vartheta_1 = 1 + \vartheta_5 \]
\[\vartheta_5\vartheta_2 = 1 + \vartheta_1 \]

\begin{figure}[h!]   
  \sidecaption
  \centerline{\includegraphics[scale= 0.5]{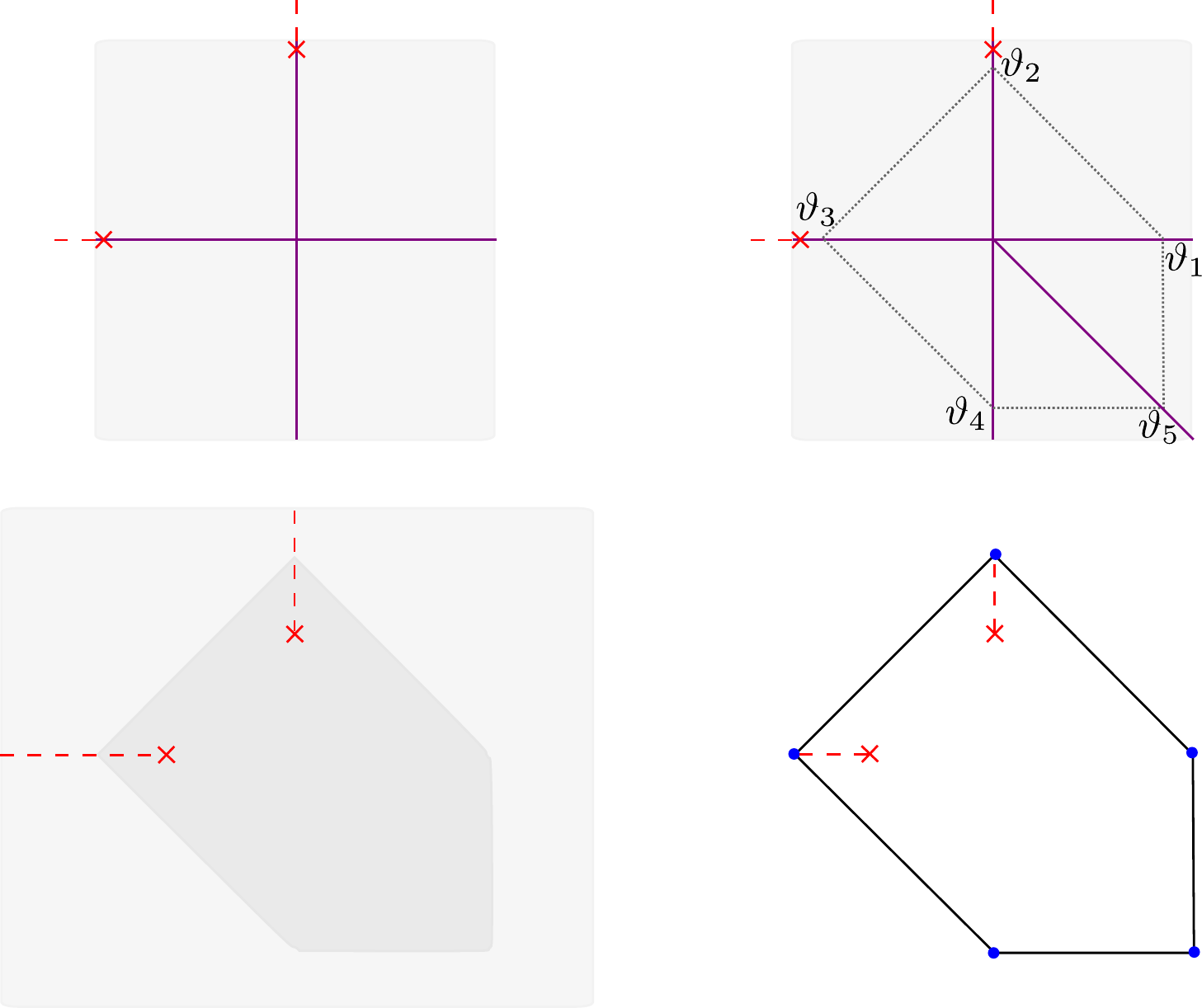}}
\caption{Scattering diagram and ATF for the $A_2$ cluster variety.}
\label{fig:dp5_00} 
\end{figure}

We see that this variety is obtained by gluing five algebraic tori
$(\C^*)^2$, with coordinates $(\vartheta_i, \vartheta_{i+1})$ [indices taken $\mod$ 5], according to the above cluster relations. We saw in more details in Section \ref{sec:ag} 
that these relations are encoded by the data of a scattering diagram, as illustrated in the top-right picture of Figure \ref{fig:dp5_00}.

Let's start with the data of an ATF describing a symplectic manifold $X$, with 2 nodal fibres, whose monodromies are encoded by cuts pointing away from the nodes in the directions $(0,1)$ and $(-1,0)$, respectively, as illustrated in the bottom-left picture of Figure \ref{fig:dp5_00}. We now think think of this as
endowed with the completed infinite volume symplectic form, so the base diagram covers the whole $\R^2$. As indicated in the previous Section, to build complex charts on this space, we add one wall for each node, represented by a line in the invariant
direction of the monodromy. [These walls represent dual fibres in the mirror $X^\vee$, bounding Maslov index 0 disks with respect to a limit complex structure $j_\infty$.] We call the chamber containing the nodes the main chamber, and we associate to 
it a complex torus $(\C^*)^2$, with coordinates $(\vartheta_2, \vartheta_{3})$, and associated with the corresponding wall 
is a gluing function of the form $(1+\vartheta_i)$. One can check that changing coordinates around these 4 walls in a full circle, does not give you identity on the  $(\vartheta_2, \vartheta_{3})$ algebraic torus. To correct for that one needs to add an extra \emph{slab}, in this case corresponding to a ray in direction $(1,-1)$, and a corresponding transition function giving you now $5$ chambers, each corresponding to an algebraic torus, 
as illustrated in the top-right picture of Figure \ref{fig:dp5_00}. This collection of walls and slabs is called the scattering diagram \cite{GS} [recall the details in Section \ref{sec:ag}]. This scattering diagram describe the relations of the $A_2$
cluster variety given in the beginning of this Section. 

We can compactify this $A_2$ cluster variety by homogenizing its defining equations, as $ \vartheta_1\vartheta_3 = \vartheta_0^2 + \vartheta_2\vartheta_0, \dots, \vartheta_5\vartheta_2 = \vartheta_0^2 + \vartheta_1\vartheta_0$. As mentioned in 
\cite[Example~8.40]{ghkk}, this gives a del Pezzo surface of degree 5 in $\CP^5$. Intersecting the hyperplane $\vartheta_0 = 0$, 
we see a chained loop of 5 divisors. Symplectically, it is natural to endow the del Pezzo surface with the monotone symplectic form
given by restricting the Fubini-Study form of $\CP^5$. The complement of the five above mentioned divisors can be seen as a 
(Weinstein) subdomain of $X$, whose completion give $X$. Indeed, there is an ATF on the degree 5 del Pezzo surface, as illustrated
in the bottom-right diagram of Figure \ref{fig:dp5_00}. [We can obtain this ATF by performing a monotone blowup in a corner of 
\cite[diagram $(A_3)$ of Figure~16]{Vi16a}.] The chained loop of 5 divisors is identified with the boundary of this ATF, and 
the complement of them is a subdomain of $X$ as illustrated by the bottom-left diagram of Figure \ref{fig:dp5_00}.


\subsection{Compactifications of Cluster varieties} \label{sec:symp_comp}

We saw in the previous section how to relate the data representing an almost-toric fibration in a open symplectic manifold 
with a set of initial walls, out of which Gross-Siebert \cite{GS} explains how to complete to a scattering diagram that provides
this manifold with complex charts given by gluing algebraic tori $(\C^*)^2$ along the walls.


\begin{figure}[h!]   
  \begin{center} 
 \centerline{\includegraphics[scale= 0.33]{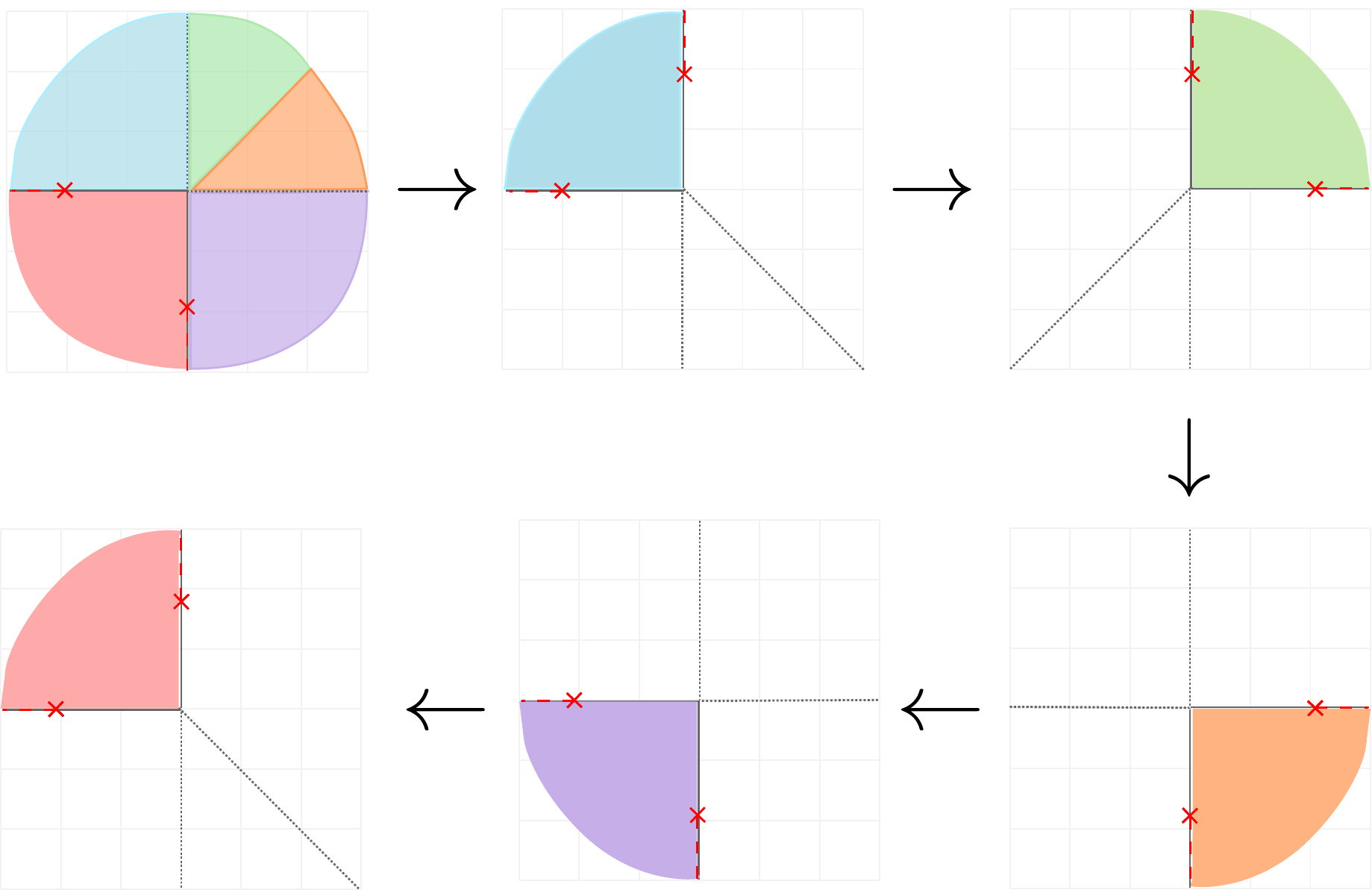}}
 \vspace{0.8cm}
  \centerline{\includegraphics[scale= 0.33]{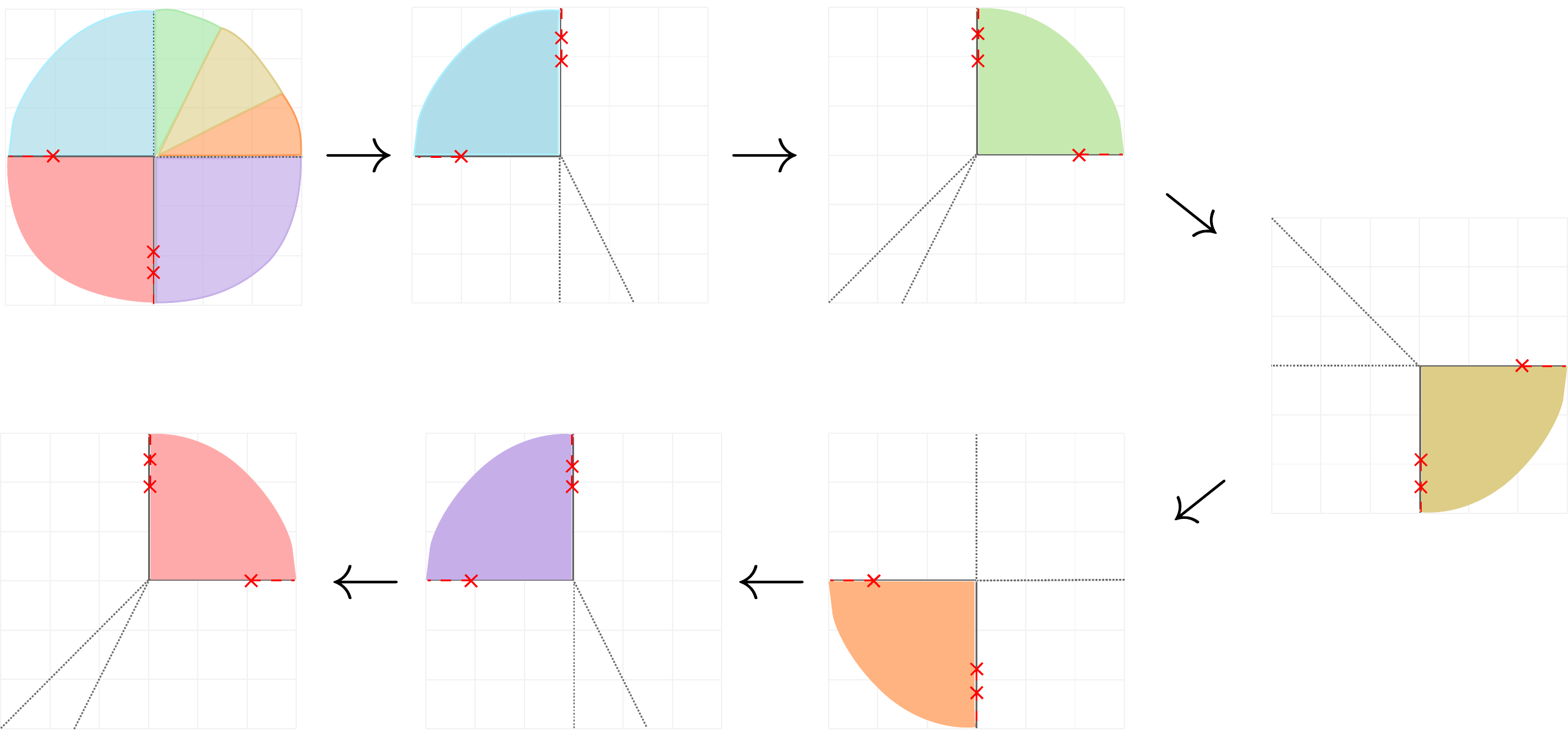}}
   \vspace{0.8cm}
  \centerline{\includegraphics[scale= 0.33]{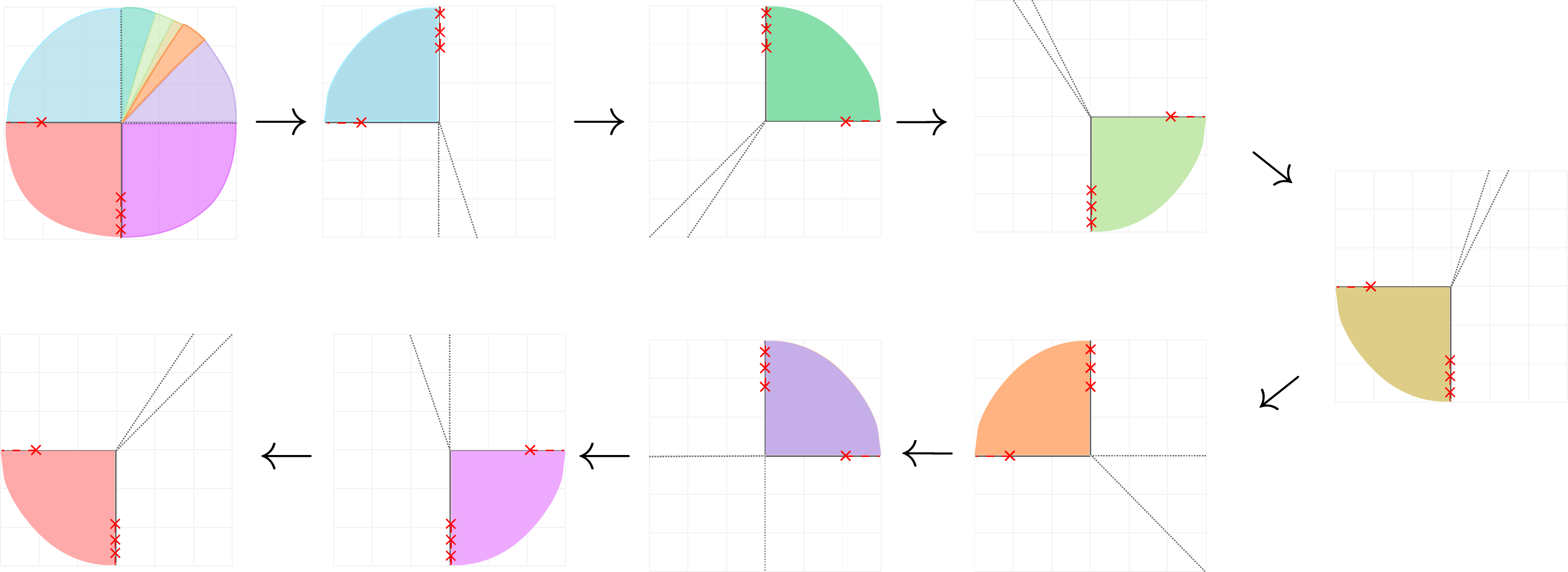}}
\caption{Cluster charts via symplectic mutations on affine cluster varieties}
\label{fig:ClusterCharts} 
\end{center} 
\end{figure}

As mentioned in Section \ref{sec:ClusterVar}, we can construct cluster varieties out of this data, and we will focus on the 
varieties of finite type $A_2$, $B_2$, $G_2$. These are open exact almost toric manifolds, built out of the scattering diagram
with initial data given by two orthogonal walls, one of them associated to one node and the other with 
one, two and three nodes, respectively, as indicated in Figure \ref{fig:ClusterCharts}. Recall we call the chart containing all nodes the main chart. One sees that symplectic mutation can be associated to changing the main chart, as illustrated in 
Figure \ref{fig:ClusterCharts}. In other words, without the prior knowledge, the 
scattering diagram can be recovered by keeping track of the ``main charts" as we apply the corresponding mutations, as illustrated in
Figure \ref{fig:ClusterCharts}. [This is not the case when the scattering diagram has a dense regions of slabs.
For instance, when considering the scattering diagram associated to the mirror of the complement of an elliptic curve in $\C \mathbb{P}^2$.
Note that this case is considered in \cite{lin2020enumerative}.]

In this section, we are interested in understanding compactifications of these cluster varieties from the symplectic perspective.
Compact symplectic manifolds have finite volume, hence we will consider as $X$ a subdomain, whose completion is the manifold described by the ATF with base diagram covering the whole $\R^2$. We will consider equivalent the subdomains with same completion.

All the symplectic compactifications considered here are symplectic del Pezzos, in the sense that they are endowed with a monotone symplectic form, which is unique up to scaling and symplectomorphisms \cite{MD90,LiLiu95,OhtaOno96,OhtaOno97}. This ensures the existence of a monotone fibre, that can be detected by the intersection point of the lines in the diagrams that go through the nodes and are in the direction of the cuts. A symplectomorphism class invariant of these monotone fibres (the star-shape) is shown \cite{ShToVi18} to be given by the interior of the polytope seen in $H^1(T^2,\R) \cong \R^2$, as we forget the nodes and cuts. So there is a symplectomorphism identifying two monotone fibres of an ATF, if and only if, the associated polytopes are related under $\mathrm{SL}(2;\Z)$. If there exists such ambient symplectomorphism, we say the Lagrangians are symplectomorphic. 

In the definition of mutation of ATFs on del Pezzo surface \cite{Vi16a}, 
besides mutating the polytope by changing the direction of the 
cut, it is required that we slide the cut through the monotone fibre. In that sense, we say that the corresponding monotone fibres are related by mutation. We can then form a graph with 
vertices representing symplectomorphism class of a Lagrangian
and edges represented these Lagrangians being related by mutation.
One aspect we can extract from the cyclic behaviour of these finite cluster varieties is the existence of cycles in the above mentioned graph.
This behaviour does not appear in mutations of monotone almost toric fibres in $\CP^2$, 
and conjecturally in $\PxP$.

\begin{figure}[h!]   
  \sidecaption
  \centerline{\includegraphics[scale= 0.4]{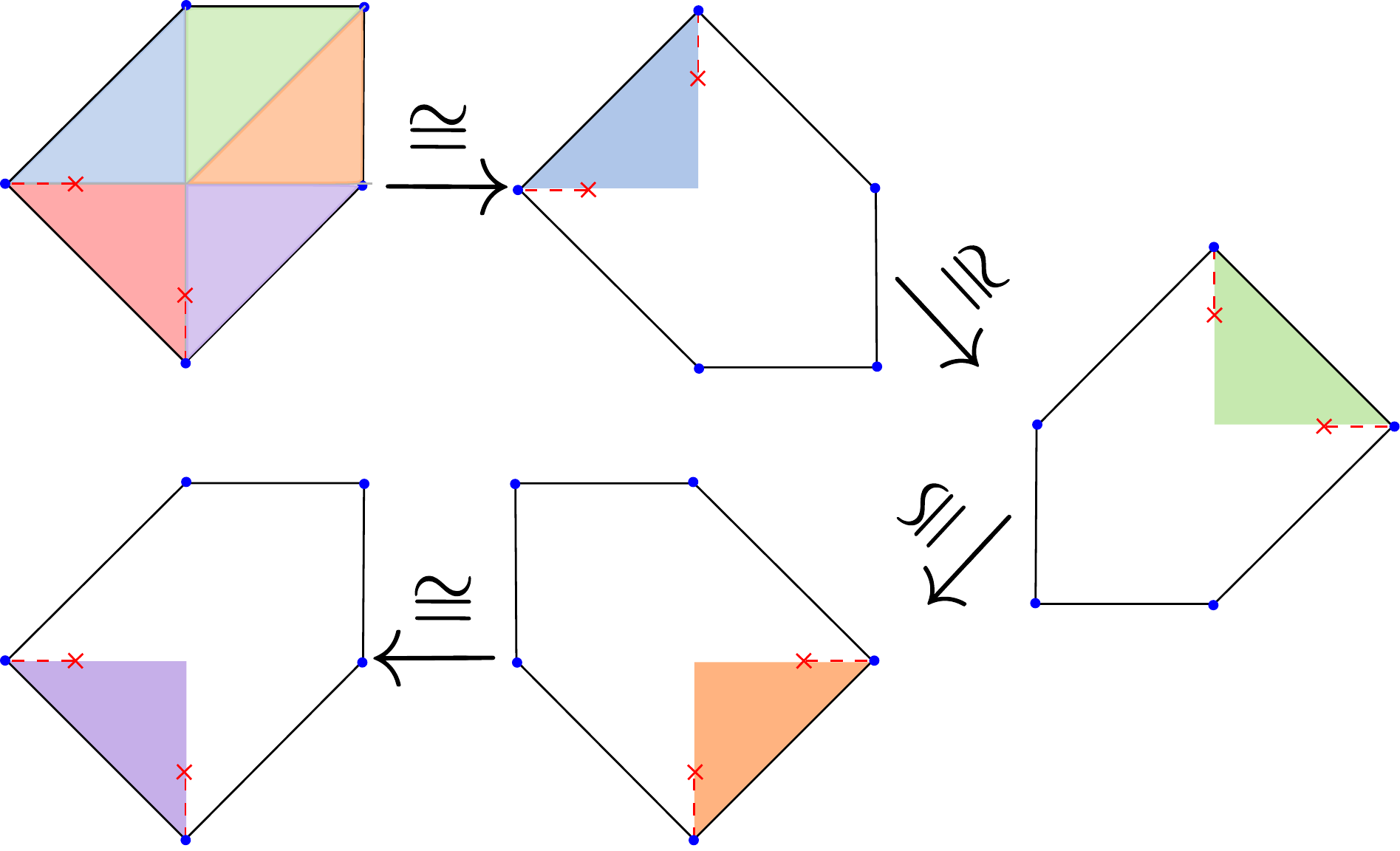}}
\caption{Mutations on degree 5 del Pezzo -- 1 torus}
\label{fig:dp5_2_1} 
\end{figure}

\begin{figure}[h!]   
  \sidecaption
   \centerline{\includegraphics[scale= 0.4]{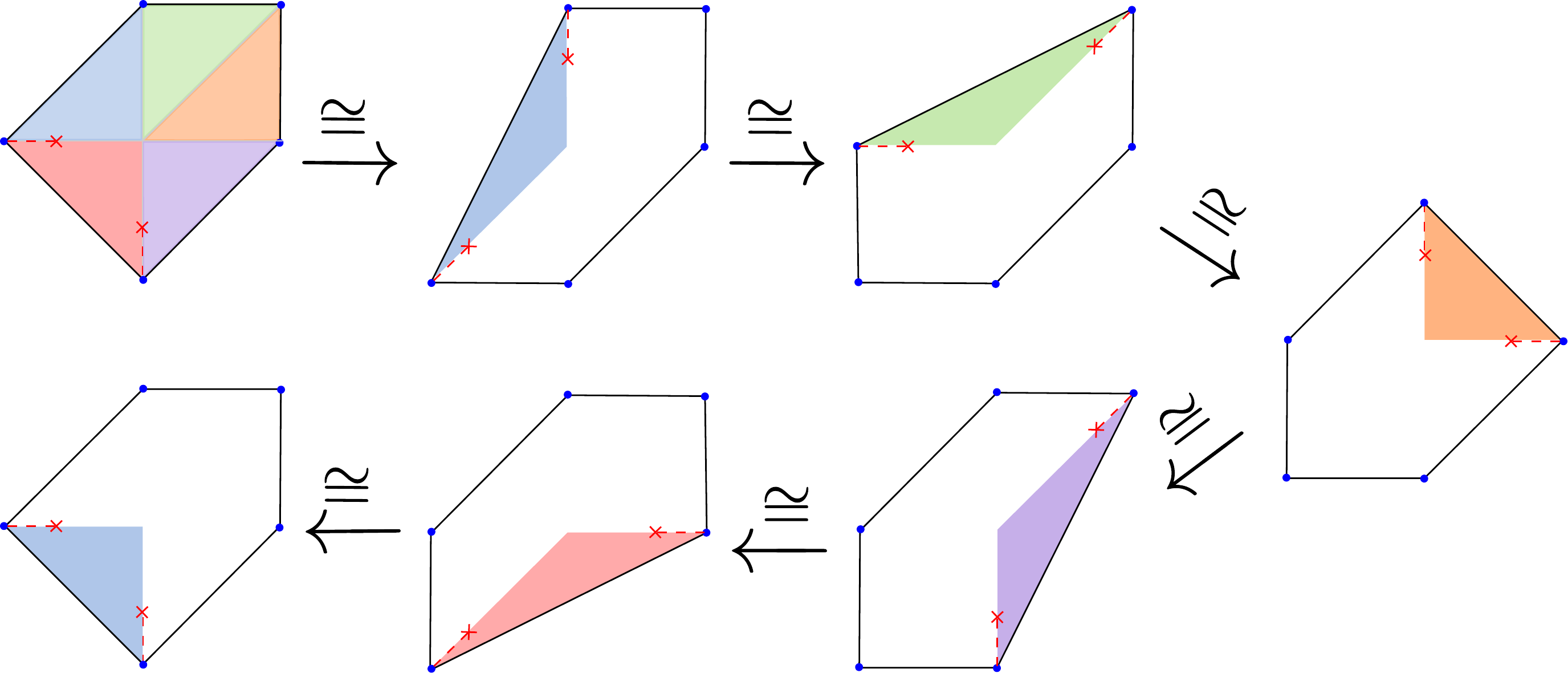}}
\caption{Mutations on degree 5 del Pezzo -- 1 torus}
\label{fig:dp5_2_2} 
\end{figure}

Let us start looking at the example from Section \ref{sec:A2} (\cite[Example~8.40]{ghkk}, \cite[Figure~1]{cpt}), a compactification of the $A_2$ cluster variety to the degree 5 del Pezzo, by adding a chain of 5 divisors, whose union represents the anti-canonical class. This compactification and its mutations are illustrated in Figures \ref{fig:dp5_2_1}, \ref{fig:dp5_2_2}. Note that we get the same pattern as in Figures \ref{fig:Xpolymutate_1}, \ref{fig:whatever}, where we get
back to the same picture after, respectively, 4 and 6 cycles, depending on the pattern of mutation. This is misleading, as ATFs, the mutations
should not depend on which half-space is fixed, and which you decide to shear. In fact, in this example, all diagrams are ($SL(2;\Z)$) equivalent,
which in particular implies that the monotone tori in each pictures are mutually symplectomorphic. The 5-cycle pattern of the $A_2$ cluster appears
by looking at the main charts, which we have already illustrated in Figure \ref{fig:ClusterCharts}. In particular, this example does not
give us a cycle of monotone Lagrangian tori, since we quotient out the graph associated to mutations by equivalence.

We want to extend a bit our notion of compactification. We will say that a (Weinstein) domain $X$ compactifies to $\bY$, if we have
$X \subset Y \subset \bY$, with $X$ a sub-domain of $Y$ and $Y = \bY \setminus \bigcup_i D_i$, for symplectic divisors $D_i$.
This will be used by us to identify our domains of interest, described in Figure \ref{fig:ClusterCharts}, appearing as open
pieces of ATFs in del Pezzo surfaces, where we do not include all the nodes. See for instance, Figures \ref{fig:dp5_1}, \ref{fig:dp4_2_1}, \ref{fig:dp3_1}, \ref{fig:dp4_1}. In these cases, our domain of interest $X$ is not the complement $Y$ of
the symplectic divisors projecting over the boundary of the ATF, but rather a subdomain of $Y$. The nodes not contained in the ATF describing $X$ will be considered frozen (not used to mutate), and depicted as a blue $\times$.

\begin{figure}[h!]   
  \begin{center} 
  \centerline{\includegraphics[scale= 0.6]{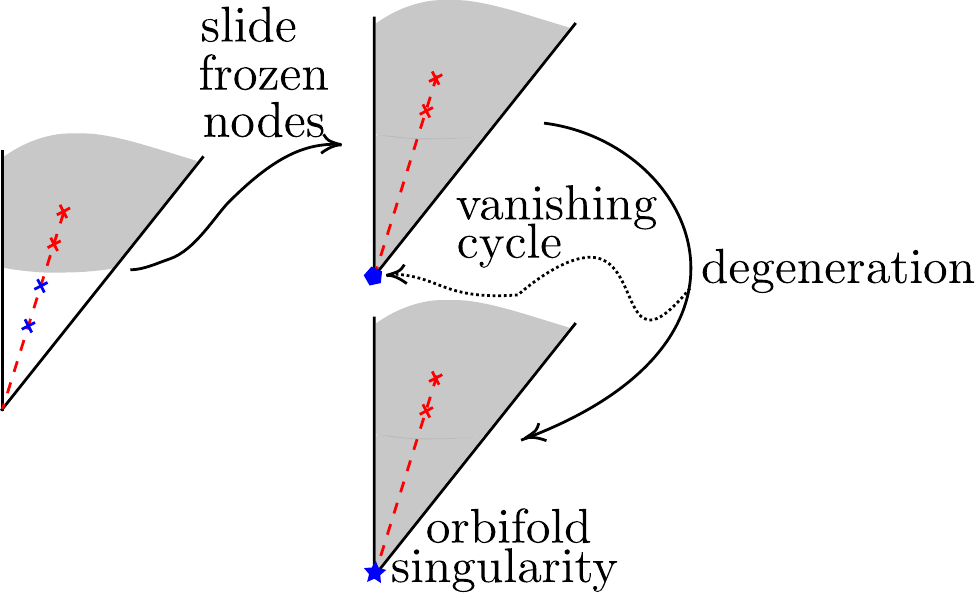}}
\caption{Moving frozen nodes to the boundary, is equivalent to have over the vertex a possibly singular Lagrangian representing a vanishing cycle of a degeneration to a toric orbifold singularity. This vanishing cycle is represented by a pentagon over the vertex, and the orbifold singularity by a star in the above diagrams.
Up to equivalence, the shaded domain can be viewed either as a subdomain of the complement of the boundary divisors in the left-picture, or the complement of the singular divisors in the orbifold diagram. }
\label{fig:Singul} 
\end{center} 
\end{figure}

An alternative way of thinking is to disregard the frozen nodes. The total manifold $\tY$ becomes singular, and a non-smooth compactification of $X$, given by adding the boundary divisor. The singularities are orbifold $T$-singularities \cite{AkKa16} at each vertex, that were previously associated with the frozen nodes. Our original smooth manifold, that included the frozen nodes, is a smoothing of this orbifold. There is a 
continuous way of relating the Lagrangian fibrations on the orbifold with the ATF on the smoothing. We like to interpret it as a two step process, which is locally illustrated in Figure \ref{fig:Singul}. The first, we keep the symplectic form on $\bY$, and consider almost-toric fibrations $ATF_t$, $t \in [0,1)$, so that in the limit $t \to 1$ the blue nodes slide all the way to a limit vertex at the boundary, and we are left with a singular Lagrangian fibration $SLF_1$, on $\bY$, such that over the limit vertices live a possibly singular Lagrangian. 
The Lagrangian over each limit vertex can be recognised in each $ATF_t$, $t < 1$ living over the associated cut from the boundary of the ATF
up to the farthest blue node and intersecting each fibre over the cut in a collapsing cycle for the corresponding nodes. 
This limit can be made rigorous but is beyond the scope of this article. The second step is
to consider a degeneration from $\bY$ to $\tY$, and a family of singular Lagrangian fibrations $SLF_s$, in the fibers corresponding
to $s \in (0,1]$, where in the limit $s \to 0$ the singular Lagrangian over each vertex collapses to the corresponding orbifold singularity. 
The Lagrangian fibrations $SLF_s$ are identified under symplectic parallel transport, so the singular Lagrangian over each vertex degenerating to an orbifold singularity is precisely the vanishing cycle of that orbifold singularity.  
In the examples presented here, the singularities associated with the frozen nodes will always be of $A_n$ type, and hence the corresponding Lagrangian a chain of $n-1$ spheres.



\begin{figure}[h!]   
  \begin{center} 
  \centerline{\includegraphics[scale= 0.4]{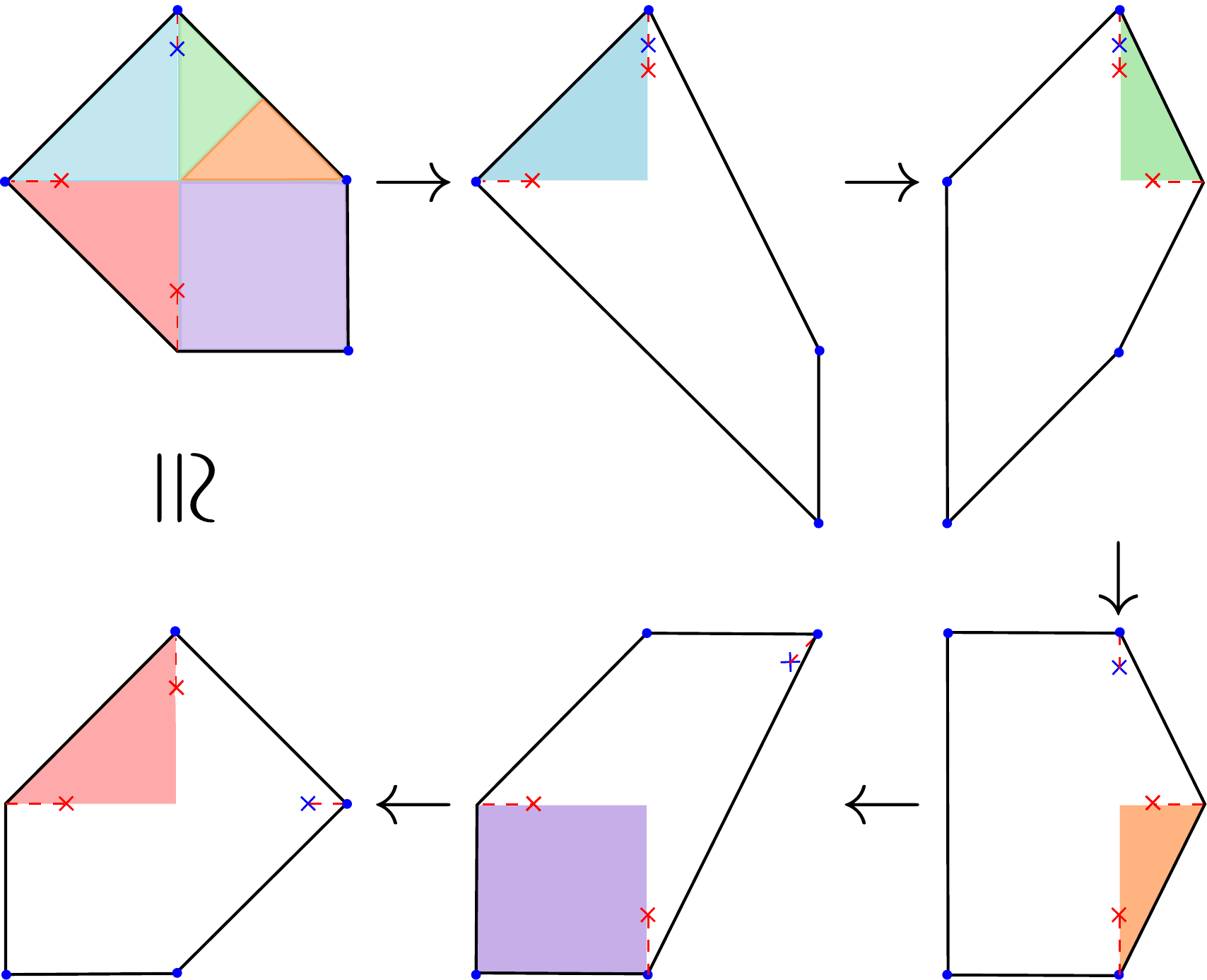}}
\caption{Mutations on degree 5 del Pezzo -- 5 tori}
\label{fig:dp5_1} 
\end{center} 
\end{figure}

Let's turn our attention now to Figure \ref{fig:dp5_1}, where we realize the degree 5 del Pezzo surface $\bY$ as a compactification of the $A_2$ cluster variety $X$, in a different way. We perform a nodal trade in a vertex at the bottom of the 
second diagram in Figure \ref{fig:dp5_2_1}, and we freeze the top node. Now the boundary divisor represents 4 symplectic spheres, 
and $X$ is a subdomain of the complement of these divisors. In this case, the mutation cycle induced by the nodal singularities in $X$ does provides us with a 5-cycle of distinct monotone Lagrangian tori. Recall that monotone fibres of non-$SL(2,\Z)$ related diagrams 
are distinct \cite{ShToVi18}. 

\begin{remark}
Disregarding the frozen node creates a double point singularity. In contrast with \cite[Example~8.40]{ghkk}, this is the 
same as considering $X_1 = 0$ in their setting. Now, consider the cycle of 5 divisors given by $\vartheta_0 = 0$. We claim that if one smooths one node of this chain (represented by our nodal trade), and then delete the resulting chain of 4 divisors in this orbifold, one recovers $X$, the $A_2$ cluster variety.
\end{remark}

\begin{figure}[h!]   
  \begin{center} 
  \centerline{\includegraphics[scale= 0.27]{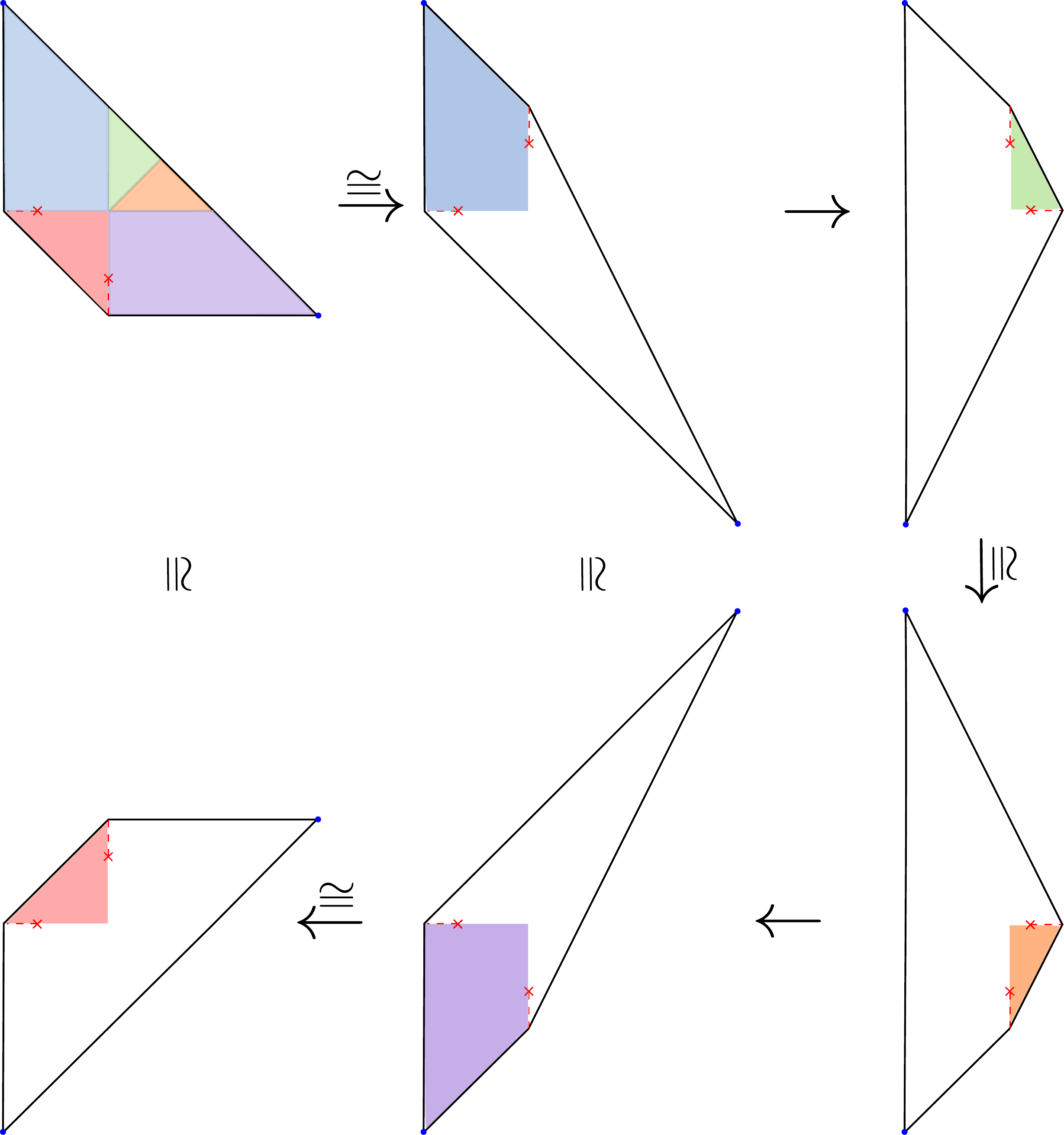}}
\caption{Mutations on degree 8 del Pezzo -- 2 tori}
\label{fig:dp8_1} 
\end{center} 
\end{figure}

We can see that there is a simpler compactification of the $A_2$ cluster variety $X$ by $\bY = \BlI$, a degree 8 del Pezzo, as illustrated in Figure \ref{fig:dp8_1} [which is the same obtained in Figure \ref{fig:monopoly}]. Here, $X$ is the complement of two divisors in classes $H$ and $2H - E$, where $H$ is the class of the line, and 
$E$ is the exceptional class. Note that we do not get a cycle of monotone Lagrangian tori, though not all tori are equivalent, we only see 2 tori in the whole cycle, which gives us only an edge on the unoriented graph of mutations of monotone tori, modulo equivalence. 


\begin{figure}[h!]   
  \begin{center} 

    \centerline{\includegraphics[scale= 0.4]{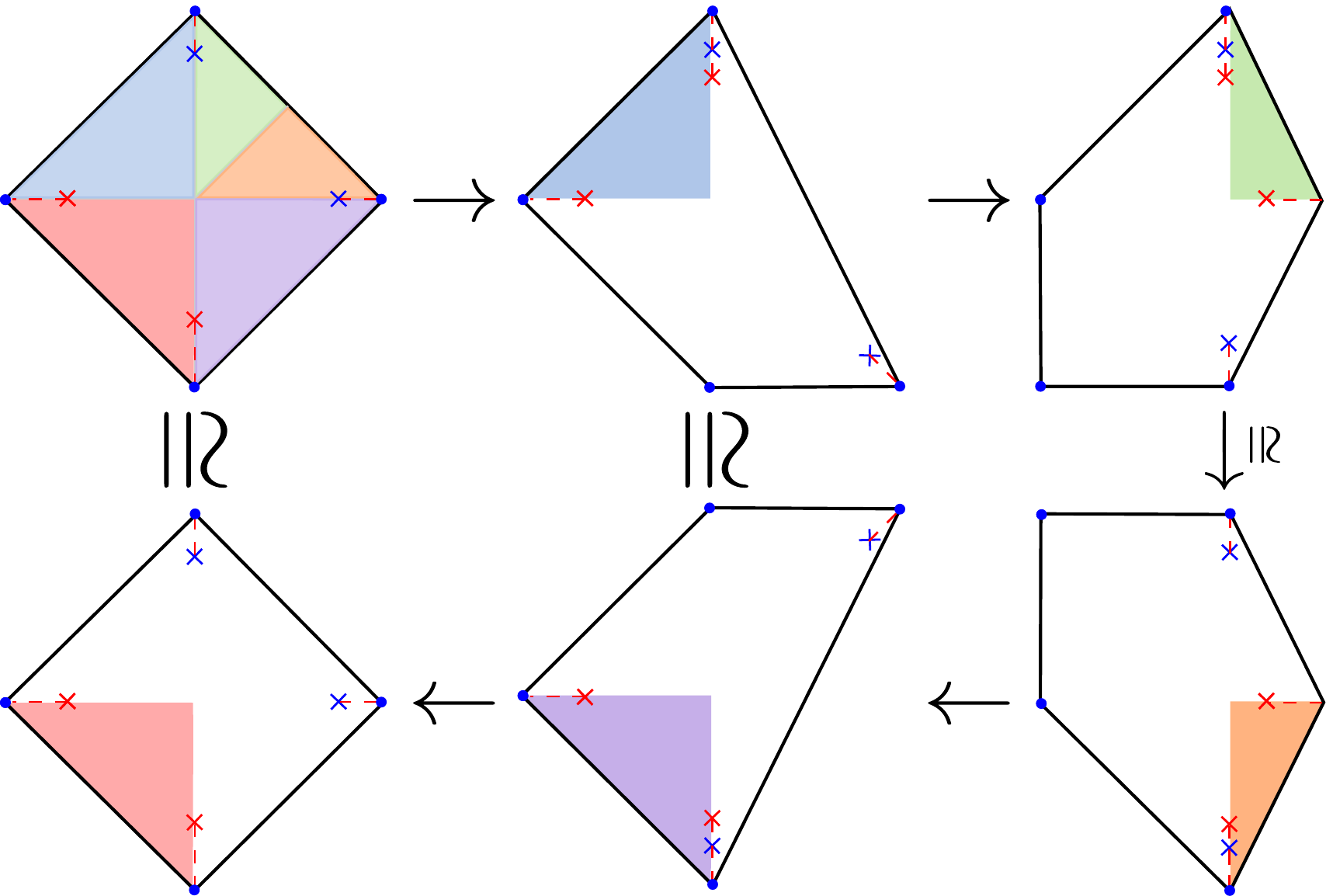}}
\caption{Mutations on degree 4 del Pezzo -- 3 tori}
\label{fig:dp4_2_1} 
\end{center} 
\end{figure}

Clearly, performing a blowup on one of the divisors of $\bY$ gives us another compactification of $X$. The top left diagram 
of Figure \ref{fig:dp4_2_1} corresponds to a toric blowup of monotone size [recall that in symplectic geometry, the blowups depend
on the size of a symplectic ball one chooses to delete] in the top left diagram of either Figure \ref{fig:dp5_2_1} or Figure \ref{fig:dp5_1}. Note in this 
case that the third and fourth, as well as the second and fifth, diagrams are equivalent, failing to deliver a cycle 
on the mutation graph of monotone Lagrangian tori in the degree 4 del Pezzo.

\begin{figure}[h!]   
\sidecaption
\includegraphics[scale= 0.35]{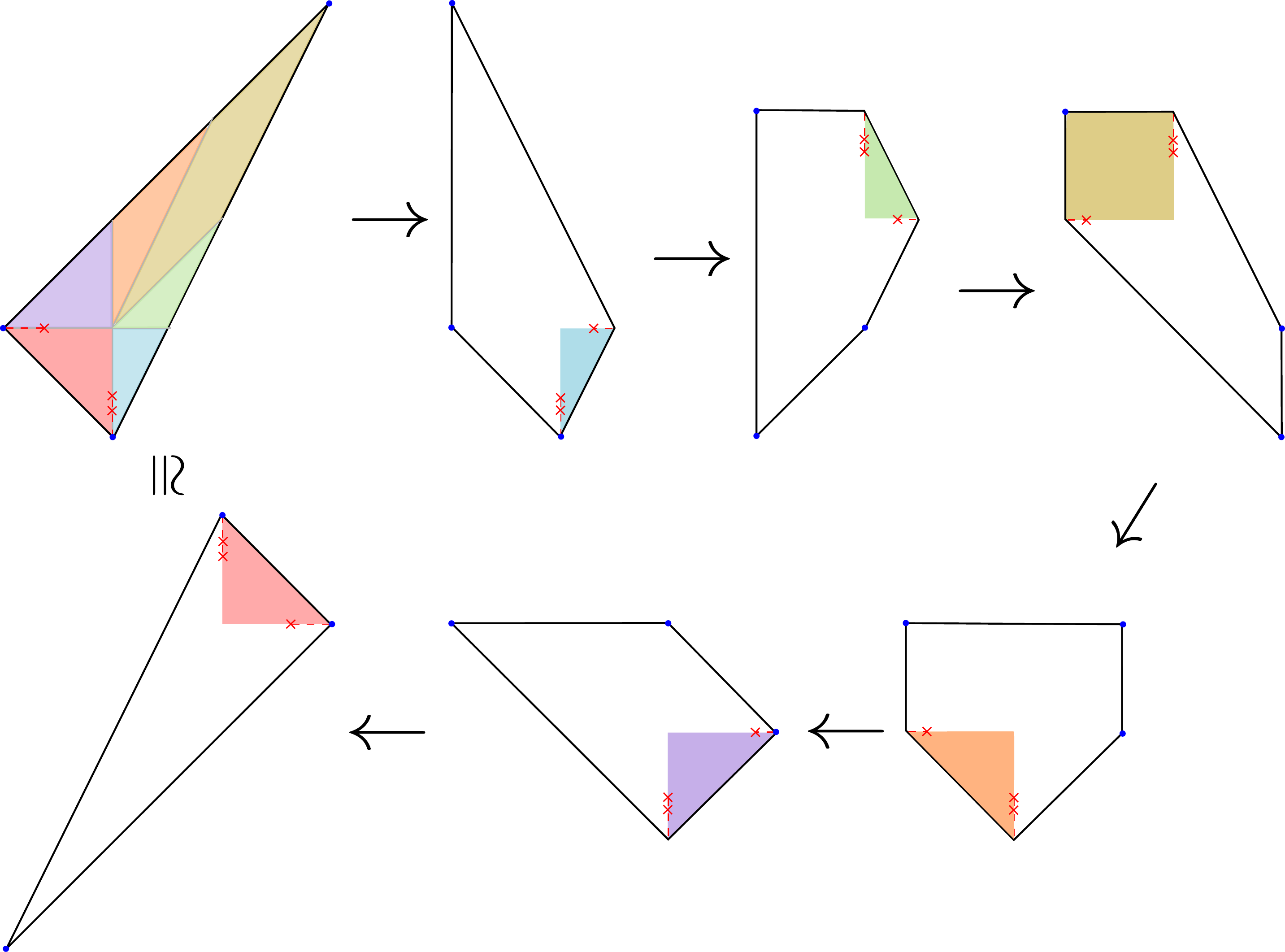}
\caption{Mutations on degree 6 del Pezzo -- 6 tori}
\label{fig:dp6_1} 

\end{figure}

\begin{figure}[h!]   
  \begin{center} 
  \centerline{\includegraphics[scale= 0.4]{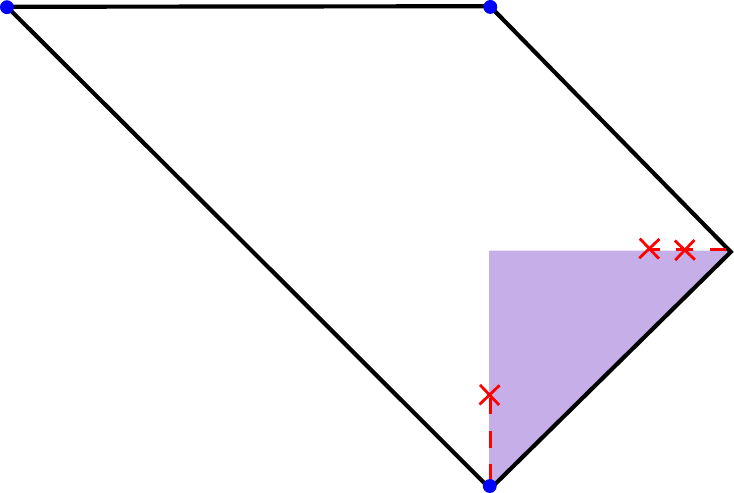}}
\caption{This diagram differ from the sixth diagram in Figure \ref{fig:dp6_1}, by one nodal trade and one inverse nodal trade. Mutations
of the displayed nodes give equivalent polytopes.}
\label{fig:dp6_Scat} 
\end{center} 
\end{figure}

\begin{figure}[h!]   
  \begin{center} 
 \centerline{\includegraphics[scale= 0.4]{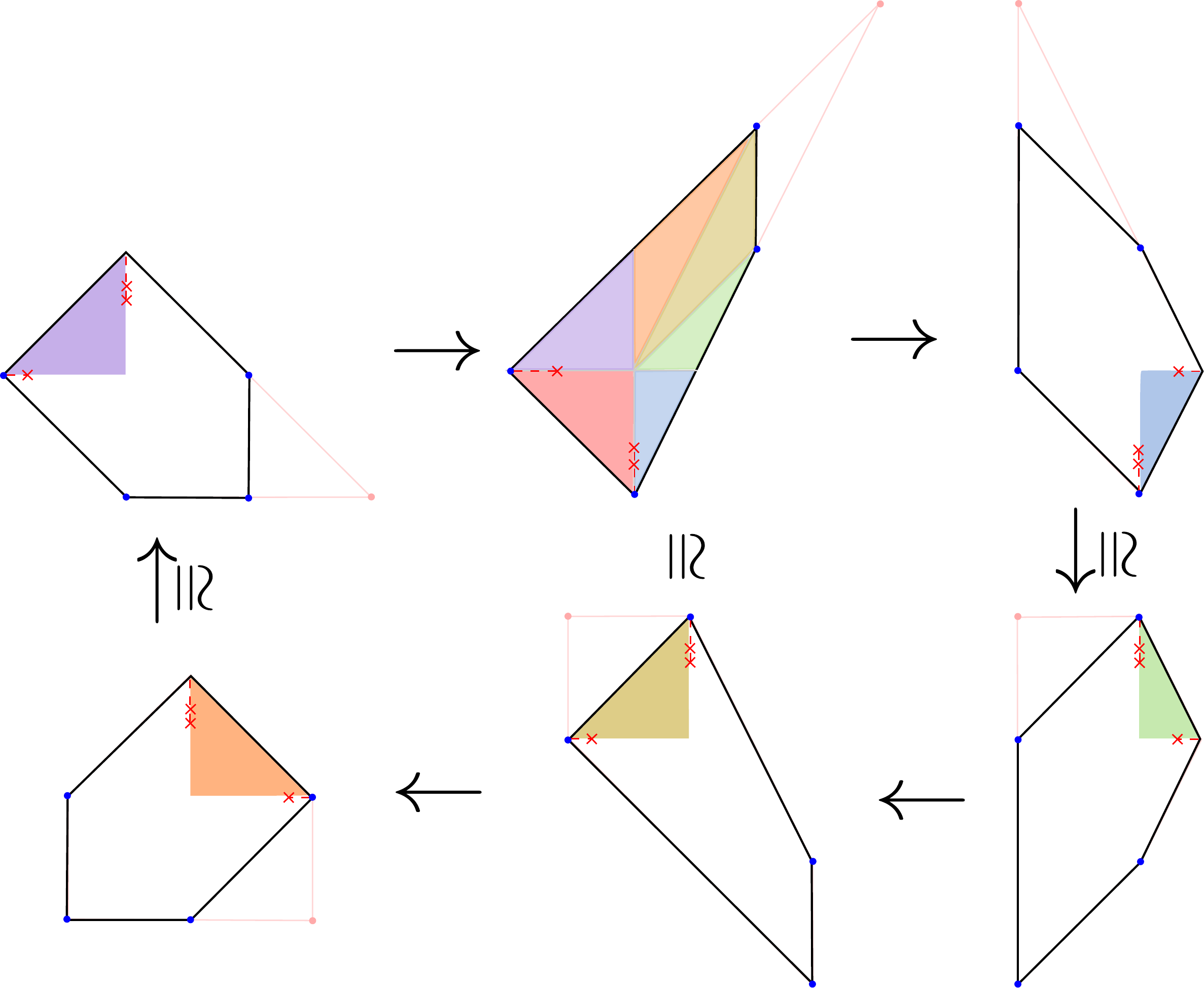}}
\caption{Mutations on degree 5 del Pezzo -- 3 tori}
\label{fig:dp5_3} 
\end{center} 
\end{figure}

Let us consider now $X$ the $B_2$ type cluster variety, with almost toric fibrations as in the series of diagrams
in the middle of Figure \ref{fig:ClusterCharts}. The first compactification $\bY$ we look at is the degree 6 del Pezzo, starting
with the ATF depicted in the top-left diagram of Figure \ref{fig:dp6_1}. [This diagram is $SL(2;\Z)$ equivalent to \cite[Diagram~$(A_5)$, Figure~16]{Vi16a} (up to nodal trades).] 
In this case $X$ is the complement of three divisors of  $\bY = \BlIII$, having symplectic areas 1, 2 and 3. We see that in this case we do get a cycle 
of size 6, with one torus corresponding to each cluster chart. 

It is interesting to notice that the sixth diagram seems to have come from 
the scattering diagram Figure \ref{fig:XB2mut_polytope}. But it is not quite the case, since that scattering diagram has a 
square function corresponding to the horizontal cut in the sixth diagram Figure \ref{fig:dp6_1}, while a simple function
corresponding to the vertical cut. This means that the natural compactification coming from that scattering diagram would be the same
del Pezzo, but represented by the diagram of Figure \ref{fig:dp6_Scat} coming from applying a nodal trade to the corner associated to the horizontal cut 
in the sixth diagram Figure \ref{fig:dp6_1}, and an inverse nodal trade on one node at the vertical cut. In particular, 
$X$ would be seen as the complement of 3 divisors in $\bY = \BlIII$, each of symplectic area 2. The reader can check that 
in this case, the mutations associated to $X$ would give equivalent monotone Lagrangian tori, analogous to the previous case depicted in 
Figures \ref{fig:dp5_2_1}, \ref{fig:dp5_2_2}.

Clearly we can also compactify the $B_2$ type cluster variety $X$ to the degree 5 del Pezzo, as the complement of four divisors as depicted
in Figure \ref{fig:dp5_3}, by simply applying a nodal trade to a diagram in Figure \ref{fig:dp5_2_1}. In Figure \ref{fig:dp5_3}, 
we depicted segments outside the diagrams to indicate that they come from applying blowups to the diagrams in 
Figure \ref{fig:dp6_1}. Curiously, it behaves similarly to the case in Figure \ref{fig:dp4_2_1}, where we have only three non-equivalent
Lagrangian tori, not providing a cycle.

\begin{figure}[h!]   
  \begin{center} 
  \centerline{\includegraphics[scale= 0.4]{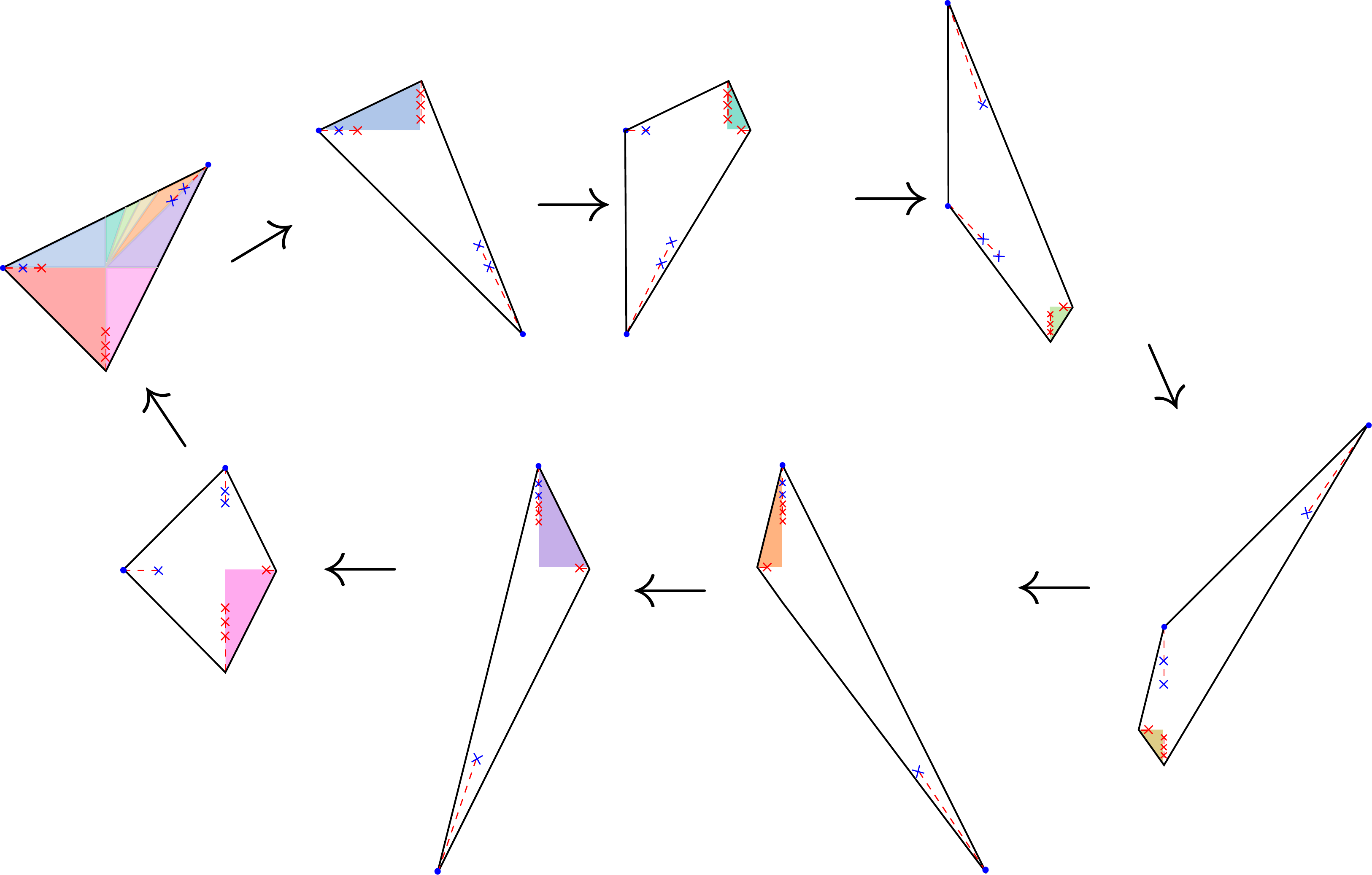}}
\caption{Mutations on degree 3 del Pezzo -- 8 tori}
\label{fig:dp3_1} 
\end{center} 
\end{figure}

\begin{figure}[h!]   
  \begin{center} 

  \centerline{\includegraphics[scale= 0.4]{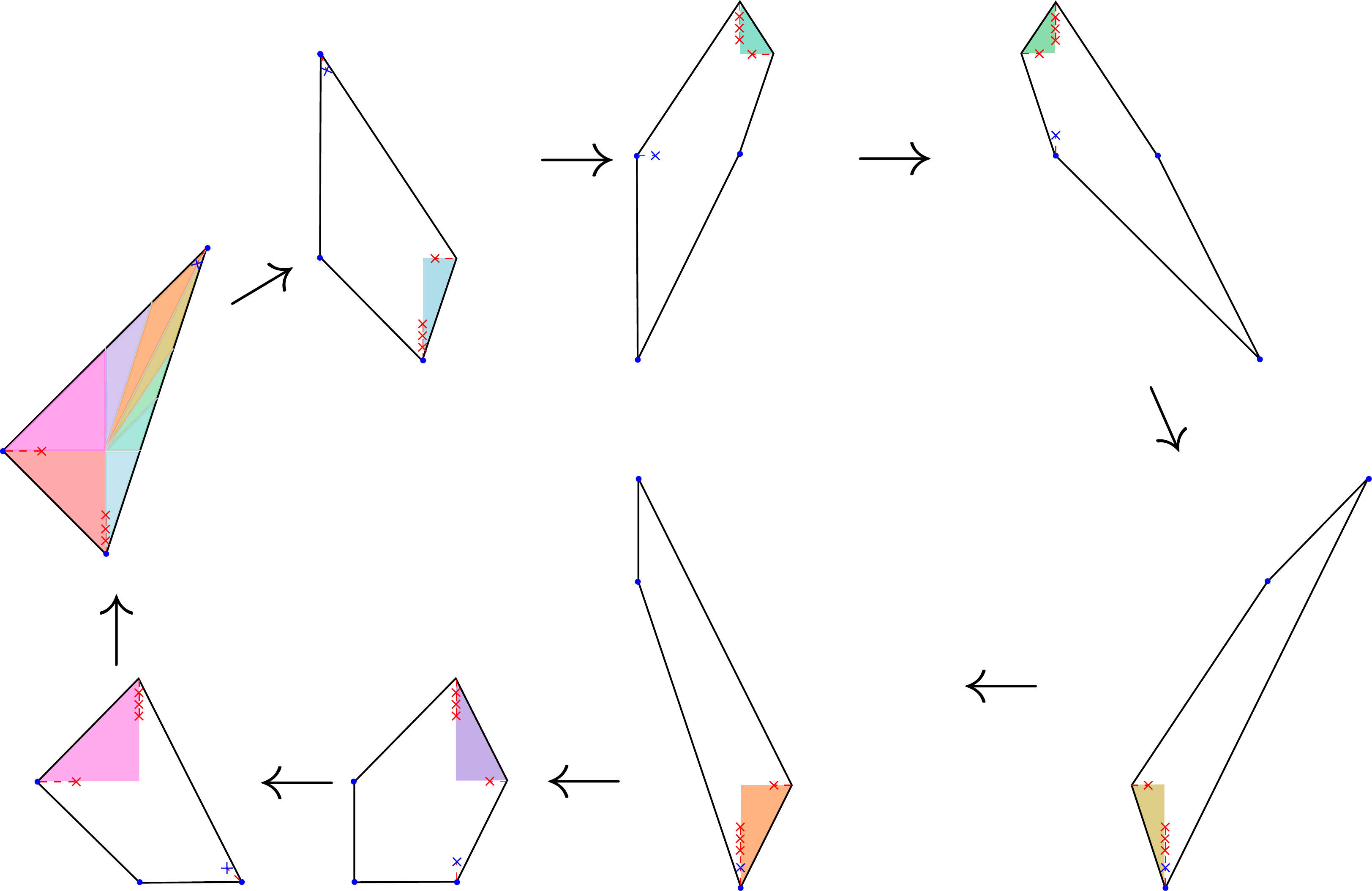}}
\caption{Mutations on degree 4 del Pezzo -- 8 tori}
\label{fig:dp4_1} 
\end{center} 
\end{figure}

We now finish by presenting two compactifications of the $G_2$ cluster variety. We name it $X$, and consider it as an almost toric variety corresponding to 
the bottom series of diagrams in Figure \ref{fig:ClusterCharts}. We start noting that the compactifications described in 
\cite{cpt}, see for instance \cite[Figure~18]{cpt}, seems to be giving partially, but not fully, smoothable orbifolds. Here we look 
to two compactifications to degree 3 and 4 del Pezzo surfaces. [Both contain frozen variables, so the reader may prefer the alternative idea of seeing $X$ compactifying to a degeneration of these surfaces.] 

We start with the top left ATBD in Figure \ref{fig:dp3_1}, which is equivalent to \cite[Diagram~$(B_2)$, Figure~19]{Vi16a}, representing an ATF of 
the cubic $\BlVI$. In this case, $X$ is a subdomain of the complement of two symplectic divisors. Sliding the 
frozen nodes to the corresponding vertex gives one Lagrangian sphere, in the horizontal cut, 
and a chain of two Lagrangian spheres in $(1,1)$-cut. This indicates that disregarding the frozen nodes corresponds to considering
an orbifold with one double-point singularity and one triple-point singularity. 
 We do get one monotone Lagrangian torus for each of the 8 cluster charts in this case.

Another compactification of $X$ is given in Figure \ref{fig:dp4_1}. The top left diagram of Figure \ref{fig:dp4_1} is equivalent (up to nodal trades) to 
\cite[Diagram~$(B_2)$, Figure~18]{Vi16a}. Here, $X$ is viewed as a subdomain of the complement of three symplectic divisors in $\bY = \BlV$. Sliding the frozen node to the vertex provides
Lagrangian sphere, or equivalently, disregarding the node gives a double-point orbifold singularity at the vertex.
As before, we get one monotone Lagrangian torus for each cluster chart.




\bigskip
%
%

\bibliographystyle{spmpsci.bst}
\bibliography{biblio} 
\end{document}